\documentclass[11pt]{amsart}
\usepackage{amsfonts}
\usepackage{amsfonts}
\usepackage[arrow,matrix]{xy}
\usepackage{amsmath,amssymb,bbm, amscd, amsthm}%,mathrsfs,hyperref}
\usepackage{longtable}
\theoremstyle{plain}
\newtheorem{thm}{Theorem}[section]
\newtheorem{cor}[thm]{Corollary}
\newtheorem{lem}[thm]{Lemma}
\newtheorem{prop}[thm]{Proposition}
\theoremstyle{definition}
\newtheorem{defn}[thm]{Definition}
\newtheorem{eg}[thm]{Example}

\newtheorem{rem}[thm]{Remark}
\def\al{\alpha}
\def\bt{\beta}
\def\dt{\delta}
\def\gm{\gamma}

\def\sg{\sigma}
\def\tt{\theta}

\def\vph{\varphi}
\def\lmd{\lambda}
\def\eps{\epsilon}
\def\vps{\varepsilon}
\def\om{\omega}

\def\bgm{\Gamma}
\def\bdt{\Delta}

\def\mc{\mathcal}

\def\ot{\otimes}
\def\op{\oplus}
\def\se{\leqslant}
\def\le{\geqslant}
\def\lan{\langle}
\def\ran{\rangle}

\def\ol{\overline}
\def\ra{\rightarrow}
\def\lra{\longrightarrow}
\def\la{\leftarrow}
\def\xra{\xrightarrow}
\def\ced{\centerdot}

\def\lap{\leftharpoonup}

\def\Hom{\operatorname {Hom}}

\def\RHom{\operatorname {RHom}}
\def\Ext{\operatorname {Ext}}

\def\ad{\operatorname {ad}}
\def\dim{\operatorname {dim}}
\def\id{\operatorname {id}}

\def\H{\operatorname {H}}
\def\HH{\operatorname {HH}}
\def\injdim{\operatorname {injdim}}

\def\hdet{\operatorname {hdet}}

\def\Cleft{\operatorname {Cleft}}

\def\t{\text}

\def\it{\textit}

\def\kk{\mathbbm{k}}
\def\ZZ{\mathbb{Z}}
\def\NN{\mathbb{N}}
\def\gg{\mathfrak{g}}

\begin{document}
\title[Cleft extensions of Koszul twisted Calabi-Yau  algebras]{\bf  Cleft extensions of Koszul twisted Calabi-Yau  algebras}

\author{Xiaolan YU}
\address {Xiaolan YU\newline Hangzhou Normal University, Hangzhou, Zhejiang 310036, China}
%\newline Department WNI, University of Hasselt, Universitaire Campus, 3590
\email{xlyu@hznu.edu.cn}

\author{Fred Van Oystaeyen}
\address{F. Van Oystaeyen\newline\indent Department of Mathematics and Computer
Science, University of Antwerp, Middelheimlaan 1, B-2020 Antwerp,
Belgium} \email{fred.vanoystaeyen@uantwerpen.be}

\author{Yinhuo ZHANG}
\address {Yinhuo ZHANG\newline Department WNI, University of Hasselt, Universitaire Campus, 3590 Diepeenbeek,Belgium } \email{yinhuo.zhang@uhasselt.be}

%\thanks{}
\date{}

\begin{abstract}
Let $H$ be a twisted Calabi-Yau (CY) algebra and $\sg$ a 2-cocycle on $H$. Let $A$ be an $N$-Koszul twisted CY algebra such that $A$ is a graded $H^\sg$-module algebra.  We show that the cleft extension $A\#_\sg H$ is also a twisted CY algebra. This result has two consequences. Firstly, the smash product of an $N$-Koszul twisted CY algebra with a twisted CY Hopf algebra is still a twisted CY algebra. Secondly, the cleft objects of a twisted CY Hopf algebra are all twisted CY algebras. As an application of this property, we determine which cleft objects of $U(\mc{D},\lmd)$, a class of pointed Hopf algebras introduced by Andruskiewitsch and Schneider, are  Calabi-Yau  algebras.
\end{abstract}


\maketitle

\section*{Introduction}
We work over a fix a field $\kk$.
Without otherwise stated, all
vector spaces, algebras are over $\kk$. Given a 2-cocycle $\sg$ on a Hopf algebra $H$ (Definition \ref{cocycle}), we can construct the algebras $H^\sg$ and $_\sg H$. Their products are deformed from the product of $H$ by
$$x*y=\sg(x_1,y_1)x_2y_2\sg^{-1}(x_3,y_3)$$
$$x\centerdot _\sg y=\sg(x_1,y_1)x_2y_2,$$
for any $x,y\in H$ respectively. The algebra $H^\sg$ together with its original coalgebra structure form a Hopf algebra, called a cocycle deformation of $H$. On the one hand, the algebra $_\sg H$   together with the original regular coaction $_\sg H\ra {}_\sg H\ot H$ form a right $H$-cleft extension over the field $\kk$. It is called a right cleft object. On the other hand, $_\sg H$ is a left $H^\sg$-cleft object with respect to the original coalgebra $_\sg H\ra H^\sg \ot {}_\sg H$. Therefore, $_\sg H$ is an $(H^\sg,H)$-bicleft object. The Hopf algebra $H^\sg$ is characterized as the Hopf algebra $L$ such that $_\sg H$ is an $(L,H)$-biGalois object (\cite{sc}).

In \cite{ma}, Masuoka studied cocycle deformations and cleft objects of a class of pointed Hopf algebras. This class of algebras includes the pointed Hopf algebras $U(\mc{D},\lmd)$ of finite Cartan type introduced by Andruskiewitsch and Schneider (\cite{as3}). The Hopf algebras $U(\mc{D},\lmd)$ consists of pointed Hopf algebras
with finite Gelfand-Kirillov dimension, which are domains with finitely generated
abelian groups of group-like elements, and generic infinitesimal braiding (\cite{aa}). By results in \cite{ma}, we know that a pointed Hopf algebra $U(D,\lambda)$ and its associated graded Hopf algebra $U(D,0)$ are cocycle deformations of each other.

The Calabi-Yau (CY for short) property of the algebras $U(\mc{D},\lmd)$ are discussed in \cite{yz}. CY algebras were introduced by Ginzburg \cite{g2} in 2006. They were studied in recent years because of their applications in algebraic geometry and mathematical physics. More general than CY algebras are so-called twisted CY algebras, which form a large class of algebras possessing the similar homological properties as the CY algebras and include CY algebras as a subclass. Associated to a twisted CY algebra, there exists a so-called Nakayama automorphism. This automorphism is unique up to an inner automorphism. A twisted CY algebra is CY if and only if its Nakayama automorphism is an inner automorphism.

For the Hopf algebra $U(\mc{D},\lmd)$, both $U(\mc{D},\lmd)$ itself and its associated graded Hopf algebra $U(\mc{D},0)$ are twisted CY algebras (\cite[Theorem 3.9]{yz}).
A more interesting phenomenon is that the CY property of $U(D,\lambda)$ is dependent only on the CY property of  $U(D,0)$. In other words,  if  $U(D,0)$ is CY, then any lifting $U(D,\lambda)$ is CY.
Note that  $U(D,\lambda)$ is a cocycle deformation of $U(D,0)$.
This raises a natural question whether  a cocycle deformation of a graded pointed (twisted) CY Hopf algebra is still a (twisted) CY algebra. For a Hopf algebra $H$ and its cocycle deformation $H^\sg$,  the algebra $_\sg H$ can be viewed as the ``connection'' between $H$ and $H^\sg $ as it defines a Morita tensor equivalence between the comodule categories over the two Hopf algebras.  To understand the relation between the twisted CY property of $H$ and that of $H^\sg$, we shall first answer the question whether $_\sg H$ is a twisted CY algebra when $H$ is.

The algebra ${}_\sg H$ can be viewed as the crossed product $\kk\#_\sg H$ (the definition of a crossed product will be reviewed in Section \ref{s1cl}).  More generally, one could ask whether the crossed product $A\#_\sg H$ will be a twisted CY algebra when both $A$ and $H$ are twisted CY algebras.  In this paper, we are able to answer the question when $A$ is a graded $N$-Koszul algebra. We note here that to form an algebra $A\#_\sg H$, it is only required that $\sg$ is an invertible map in $\Hom(H\ot H,A)$ satisfying the cocycle condition and $A$ is a twisted $H$-module. When $A$ is a graded $N$-Koszul algebra, the assumption that $\sg$ has its image in $\kk$ is necessary to make sure that the obtained crossed product $A\#_\sg H$ is still a graded algebra. In this case $\sg$ is just a 2-cocycle on $H$ and $A$ is a left graded $H^\sg$-module algebra. Here $A$ is a left graded $H^\sg$-module algebra means that $A$ is a left $H^\sg$-module algebra such that each graded piece $A_i$ is a left $H^\sg$-module.  The following theorem is our main result (see Theorem \ref{main}):

\begin{thm}\label{in}
Let $H$ be a twisted CY Hopf algebra with homological integral $\int^l_H=\kk_\xi$, where $\xi:H\ra \kk$ is an algebra homomorphism and $\sg$ a 2-cocycle on $H$. Let $A$ be a $N$-Koszul graded twisted CY algebra with Nakayama automorphism $\mu$ such that $A$ is a left graded $H^\sg$-module algebra.
 Then $A\#_\sg H$ is a twisted CY algebra with Nakayama automorphism  $\rho$ defined by $\rho(a\#h)=\mu(a)\#\hdet_{H^\sg}(h_1)(S_{\sg,1}^{-1}(S_{1,\sg}^{-1}(h_2)))\xi(h_3)$ for all $a\#h\in A\#_\sg H$.
\end{thm}
Here, $\hdet_{H^\sg}$ denotes the homological determinant of the $H^\sg$-action. The homological integral of a twisted CY Hopf algebra will be given in Section \ref{s2}. The notion $S_{\sg,\tau}$ will be recalled in Section \ref{s1c}. Examples of Theorem \ref{in} will be provided in Section \ref{s4}.

Theorem \ref{in} has two consequences. Firstly, in Theorem \ref{in}, if we let the cocycle $\sg$ be trivial, then the crossed product $A\#_\sg H$ is just the smash product $A\#H$. Therefore, we obtain the following result on smash products.

\begin{thm}\label{in1}
Let $H$ be a twisted CY Hopf algebra with homological integral $\int^l_H=\kk_\xi$, where $\xi:H\ra \kk$ is an algebra homomorphism and $A$ an $N$-Koszul graded twisted CY algebra with Nakayama automorphism $\mu$ such that $A$ is a left graded $H$-module algebra. Then $A\#H$ is a twisted CY algebra with Nakayama automorphism  $\rho$ defined by $\rho(a\#h)=\mu(a)\#\hdet_{H}(h_1)(S^{-2}(h_2))\xi(h_3)$, for any $a\#h\in A\#H$.
\end{thm}
This generalizes the results in \cite{liwz} and \cite{rrz}. The smash products of CY algebras has been studied quite broadly. For instance, see \cite{f}, \cite{ir},  \cite{liwz}, \cite{wz}, \cite{rrz}. The results in \cite{liwz} and \cite{rrz} are probably two of the most general results in this direction. \cite{liwz} states that  when $H$ is  an involutory Hopf CY algebra and $A$ is an $N$-Koszul CY algebra, the smash product $A\#H$ is CY if and only if the homological determinant of the $H$-action on $A$ is trivial. One of the main results in \cite{rrz} states  that the smash product $A\#H$ is a twisted CY algebra when $A$ is a graded twisted CY algebra and $H$ a finite dimensional Hopf algebra acting on $A$. The Nakayama automorphism of $A\#H$ is determined by the ones of $A$ and $H$, along with the homological determinant of the $H$-action.

Secondly, in Theorem \ref{in}, if we let the algebra $A$ be $\kk$, we obtain the following description of the twisted CY property of cleft objects.
\begin{thm}\label{in2}
Let $H$ be a twisted CY Hopf algebra with $\int^l_H={}_\xi\kk$, and $\sg$ a 2-cocycle on $H$.    Then the right cleft object ${}_\sg H$ is a twisted CY algebra with Nakayama automorphism $\mu$ defined by
$$\mu(x)=S^{-1}_{\sg,1}(S^{-1}_{1,\sg}(x_1))\xi S(x_2)$$
for any $x\in {}_\sg H$.
\end{thm}

As an application of Theorem \ref{in2}, we study the CY property of the cleft objects of the Hopf algebras $U(\mc{D},\lmd)$ in Section \ref{s3}. It turns out that all cleft objects of the algebra $U(\mc{D},\lmd)$ are twisted CY algebras. Their Nakayama automorphisms are given explicitly in Proposition \ref{cycle}. Hence we are able to characterize when a clefts object is CY. It is interesting that a cleft object of $U(\mc{D},\lmd)$ could be a  CY algebra even when $U(\mc{D},\lmd)$ itself is not. We give such an example at the end of Section \ref{s3}.

Our motivating examples are the algebras of the form $A\#_\sg\kk G$, where $A$ is a polynomial algebra, $G$ is a finite group acting on $A$, and $\sg:G\times G\ra \mathbb{C}^\times$ is a 2-cocycle on $G$. Such crossed products are of interest in geometry due to their relationship with corresponding orbifolds (for e.g., see \cite{ar}, \cite{cgw}, \cite{vw}). In Section \ref{s4}, we show that these crossed products are all twisted CY algebras.  PBW deformations of the crossed product $A\#_\sg \kk G$ are the twisted Drinfeld Hecke algebras defined in \cite{wi}. If the cocycle is trivial, then $A\#\kk G$, the skew group algebra, is just the Drinfeld Hecke algebras defined by V. Drinfeld \cite{dr}. They  have been studied by many authors, for example \cite{eg}, \cite{bb}, \cite{l1}.  Quantum Drinfeld Hecke algebras are anther generalizations of Drinfeld Hecke algebras by replacing polynomial algebras by quantum polynomial algebras \cite{ls}, \cite{nw}. More generally, Naidu defined twisted quantum Drinfeld Hecke algebras in \cite{na}.  A twisted quantum Drinfeld Hecke algebra is an algebra of the form $A\#_\sg\kk G$, where $A$ is a quantum polynomial algebra, $G$ is a finite group acting on $A$, and $\sg$ is a 2-cocycle on $G$. Twisted quantum Drinfeld Hecke algebras are generalizations of both twisted Drinfeld Hecke algebras and quantum Drinfeld Hecke algebras.
A quantum polynomial algebra is a Koszul algebra. If  PBW deformations of the algebra $A\#_\sg H$ in Theorem \ref{in} are still twisted CY algebras, then twisted quantum Drinfeld Hecke algebras will all be twisted CY algebras. We will discuss this problem in our upcoming paper.

\section{Preliminaries }\label{s1}

Throughout this paper, the unadorned  tensor $\ot$
means $\ot_\kk$ and $\Hom$ means $\Hom_\kk$.

Given an algebra $A$, we write $A^{op}$ for the
opposite algebra of $A$ and $A^e$ for the enveloping algebra $A\ot
A^{op}$. An $A$-bimodule can be identified with a left $A^e$-module or a right $A^e$-module.

For an $A$-bimodule $M$ and two algebra automorphisms $\mu$ and $\nu$, we let $^\mu M^\nu$ denote the $A$-bimodule such that $^\mu M^\nu\cong M$ as vector spaces, and the bimodule structure is given by
$$a\cdot m \cdot b=\mu(a)m\nu(b),$$
for all $a,b\in A$ and $m\in M$. If one of the automorphisms is the identity, we will omit it. It is well-known that  $A^\mu\cong {}^{\mu^{-1}} A$ as
$A$-$A$-bimodules. $A^\mu\cong A$ as $A$-$A$-bimodules if and only
if $\mu$ is an inner automorphism.

We assume that the Hopf algebras considered in this paper have bijective antipodes. For a Hopf algebra $H$, we use  Sweedler's (sumless) notation for the
comultiplication and coaction of $H$.

\subsection{Cogroupoid}\label{s1c}
\begin{defn}

a \it{cocategory} $\mc{C}$ consists of:
\begin{itemize}
\item A set of objects $\t{ob}(\mc{C})$.
\item For any $X,Y\in \t{ob}(\mc{C})$, an algebra $\mc{C}(X,Y)$.
\item For any $X,Y,Z\in \t{ob}(\mc{C})$, algebra homomorphisms
$$\bdt^Z_{XY}:\mc{C}(X,Y)\ra \mc{C}(X,Z)\ot \mc{C}(Z,Y) \t{ and } \vps_X:\mc{C}(X,X)\ra \kk$$
such that for any $X,Y,Z,T\in \t{ob}(\mc{C})$, the following diagrams commute:
$$\begin{CD}\mc{C}(X,Y)@>\bdt^Z_{X,Y}>>\mc{C}(X,Z)\ot \mc{C}(Z,Y)\\
@V\bdt^T_{X,Y}VV@V\bdt^T_{X,Z}\ot 1VV\\
\mc{C}(X,T)\ot \mc{C}(T,Y)@>1\ot \bdt^Z_{T,Y}>>\mc{C}(X,T)\ot \mc{C}(T,Z)\ot \mc{C}(Z,Y)\end{CD}$$
\xymatrix
{\mc{C}(X,Y)\ar@{=}[rd]\ar^{\bdt^Y_{X,Y}}[d]&\\
\mc{C}(X,Y)\ot\mc{C}(Y,Y)\ar^-{1\ot \vps_Y}[r]&\mc{C}(X,Y)}\hspace{2mm}
\xymatrix
{\mc{C}(X,Y)\ar@{=}[rd]\ar^{\bdt^X_{X,Y}}[d]&\\
\mc{C}(X,X)\ot\mc{C}(X,Y)\ar^-{\vps_X \ot 1 }[r]&\mc{C}(X,Y).}
\end{itemize}
\end{defn}
Thus a cocategory with one object is just a bialgebra.

A cocategory $\mc{C}$ is said to be \it{connected} if $\mc{C}(X,Y)$ is a non zero algebra for any $X,Y\in \t{ob}(\mc{C})$.

\begin{defn}A \it{cogroupoid} $\mc{C}$ consists of a cocategory $\mc{C}$ together with, for any $X,Y\in \t{ob}(\mc{C})$, linear maps
$$S_{X,Y}:\mc{C}(X,Y)\longrightarrow \mc{C}(Y,X)$$
such that for any $X,Y\in \mc{C}$, the following diagrams commute:
$$\xymatrix{\mc{C}(X,X)\ar[d]_{\bdt_{X,X}^Y}\ar[r]^-{\vps_X}&\kk\ar[r]^-u&\mc{C}(X,Y)\\
\mc{C}(X,Y)\ot\mc{C}(Y,X)\ar[rr]^{1\ot S_{Y,X}}&&\mc{C}(X,Y)\ot\mc{C}(X,Y)\ar[u]^m}$$
$$\xymatrix{\mc{C}(X,X)\ar[d]_{\bdt_{X,X}^Y}\ar[r]^-{\vps_X}&\kk\ar[r]^-u&\mc{C}(Y,X)\\
\mc{C}(X,Y)\ot\mc{C}(Y,X)\ar[rr]^{S_{X,Y}\ot 1}&&\mc{C}(Y,X)\ot\mc{C}(Y,X)\ar[u]^m}$$
\end{defn}

We refer to \cite{bi1} for basic properties of cogroupoids.

In this paper, we are mainly
concerned with the 2-cocycle cogroupoid of a Hopf algebra.

\begin{defn}\label{cocycle}Let $H$ be a Hopf algebra. A  \it{(right) 2-cocycle} on $H$ is a convolution invertible linear map $\sg:H\ot H\ra \kk$ satisfying
\begin{equation}\label{rcocyle1}\sg(h_1,k_1)\sg(h_2k_2,l)=\sg(k_1,l_1)\sg(h,k_2l_2)\end{equation}
\begin{equation}\label{rcocyle2}\sg(h,1)=\sg(1,h)=\vps(h)\end{equation} for all $h,k,l\in H$. The set of 2-cocycles on $H$ is denoted $Z^2(H)$.
\end{defn}

 The convolution inverse of $\sg$, denote $\sg^{-1}$, satisfies
\begin{equation}\label{lcocyle1}\sg^{-1}(h_1k_1,l)\sg^{-1}(h_2,k_2)=\sg^{-1}(h,k_1l_1)\sg^{-1}(k_2,l_2)\end{equation}
\begin{equation}\label{lcocyle2}\sg^{-1}(h,1)=\sg^{-1}(1,h)=\vps(h)\end{equation} for all $h,k,l\in H$. Such a convolution invertible map is called a \textit{left 2-cocycle} on $H$.
Conversely, the convolution inverse of a left 2-cocycle is just a right 2-cocycle.
%The equations (\ref{lcocyle1}) and (\ref{lcocyle2}) are actually the mirror image formulae of (\ref{rcocyle1}) and (\ref{rcocyle2}).

The set of 2-cocycles defines the 2-cocycle cogroupoid of $H$.

Let $\sg,\tau\in Z^2(H)$. The algebra $H(\sg,\tau)$ is defined to be the vector space $H$ together with the multiplication  given by
\begin{equation}\label{mulco}h\centerdot k=\sg(h_1,k_1)h_2k_2\tau^{-1}(h_3,k_3),\end{equation}
for any $h,k\in H$.

Now we recall the necessary structural maps for the 2-cocycle cogroupoid on $H$. For any $\sg,\tau, \om\in Z^2(H)$, define the following maps:
\begin{equation}\label{comulco}\begin{array}{rcl}\bdt_{\sg,\tau}^\om=\bdt:H(\sg,\tau)&\longrightarrow &H(\sg,\om)\ot H(\om,\tau)\\
h&\longmapsto &h_1\ot h_2.
 \end{array}\end{equation}
\begin{equation}\label{couco}\vps_\sg=\vps:H(\sg,\sg)\longrightarrow \kk.\end{equation}
\begin{equation}\label{antico}\begin{array}{rcl}S_{\sg,\tau}:H(\sg,\tau)&\longrightarrow & H(\tau,\sg)\\
h&\longmapsto &\sg(h_1,S(h_2))S(h_3)\tau^{-1}(S(h_4),h_5).\end{array}\end{equation}
It is routine to check that the inverse of $S_{\sg,\tau}$ is given as follows:
\begin{equation}\label{anticoinverse}\begin{array}{rcl}S^{-1}_{\sg,\tau}:H(\tau,\sg)&\longrightarrow & H(\sg,\tau)\\
h&\longmapsto &\sg^{-1}(h_5,S^{-1}(h_4))S^{-1}(h_3)\tau(S^{-1}(h_2),h_1).\end{array}\end{equation}

The \it{$2$-cocycle cogroupoid} of $H$, denoted by $\underline{H}$, is the cogroupoid defined as follows:
\begin{enumerate}
\item[(i)] $\t{ob}(\underline{H})=Z^2(H)$.
\item[(ii)] For $\sg,\tau\in Z^2(H)$, the algebra $\underline{H}(\sg,\tau)$ is the algebra $H(\sg,\tau)$ defined in (\ref{mulco}).
\item[(iii)] The structural maps $\bdt^\bullet_{\bullet,\bullet}$, $\vps_\bullet$ and $S_{\bullet,\bullet}$ are defined in (\ref{comulco}), (\ref{couco}) and (\ref{antico}) respectively.
\end{enumerate}
\cite[Lemma 3.13]{bi1} shows that the maps $\bdt^\bullet_{\bullet,\bullet}$, $\vps_\bullet$ and $S_{\bullet,\bullet}$ indeed satisfy the conditions required for a cogroupoid. It is clear that a 2-cocycle cogroupoid is connected. The following lemma follows from basis properties of cogroupoids.

\begin{lem}\label{anpo}\cite[Proposition 2.13]{bi1}Let $\underline{H}$ be the 2-cocycle cogroupoid, and let $\sg,\tau\in \t{ob}(\underline{H})$.
\begin{enumerate}
\item[(i)] $S_{\sg,\tau}:H(\sg,\tau)\ra H(\tau,\sg)^{op}$ is an algebra homomorphism.
\item[(ii)] For any $\omega\in \t{ob}(\underline{H})$ and $h\in H$, we have
$$\bdt_{\tau,\sg}^\omega(S_{\sg,\tau}(h))=S_{\om,\tau}(h_1)\ot S_{\sg,\om}(h_2).$$
\end{enumerate}

\end{lem}

The Hopf algebra $H(1,1)$ (where $1$ stands for $\vps\ot\vps$) is just the Hopf algebra $H$ itself. Let $\sg$ be a 2-cocycle. We write ${}_\sg H$ for the algebra $H(\sg,1)$.  Similarly, we write $H_{\sg^{-1}}$ for the algebra $H(1,\sg)$.  To make the presentation clear,  we let $\centerdot_\sg$ and $\centerdot_{\sg^{-1}}$  denote the multiplications in ${}_\sg H$ and $H_{\sg^{-1}}$ respectively.

The Hopf algebra $H(\sg,\sg)$ is just the \textit{cocycle deformation} $H^\sg$ of $H$ defined by Doi in \cite{do}. The comultiplication of $H^\sg$ is the same as the comultiplication of $H$. However, the multiplication and the antipode are deformed:
$$h*k=\sg(h_1,k_1)h_2k_2\sg^{-1}(h_3,k_3),$$
$$S_{\sg,\sg}(h)=\sg(h_1,S(h_2))S(h_3)\sg^{-1}(S(h_4),h_5)$$
for any $h,k\in H^\sg$. In the following,  $S_{\sg, \sg}$ is denoted by $S^\sg$ for simplicity.

\subsection{Cleft extensions}\label{s1cl}
A Hopf algebra $H$ is said to \textit{measure} an algebra $A$ if there is a $\kk$-linear map $H\ot A\ra A$, given by $h\ot a\mapsto h\cdot a$, such that $h\cdot 1=\vps(h)$ and $h\cdot(ab)=(h_1\cdot a)(h_2\cdot b)$ for all $h\in H$, $a,b\in A$.

\begin{defn}
Let $H$ be a Hopf algebra and $A$ an algebra. Assume that $H$ measures $A$ and that $\sg$ is an invertible map in $\Hom(H\ot H,A)$. The \it{crossed product} $A\#_\sg H$ of $A$ with $H$ is defined on the vector space $A\ot H$  with multiplication given by
$$(a\#h)(b\#k)=a(h_1\cdot b)\sg(h_2,k_1)\#h_3k_2$$
for all $h,k\in H$, $a,b\in A$. Here we write $a\#h$ for the tensor product $a\ot h$.
\end{defn}

The following lemma is well-known (cf. \cite[Lemma 7.1.2]{m}).
\begin{lem}\label{cross}
$A\#_\sg H$ is an associative algebra with identity element $1\#1$ if and only if the following two conditions are satisfied:
\begin{enumerate}
\item[(i)] $A$ is a twisted $H$-module. That is, $1\cdot a=a$ for all $a\in A$, and
$$h\cdot(k\cdot a)=\sg(h_1,k_1)(h_2k_2\cdot a)\sg^{-1}(h_3,k_3),$$
for all $h,k\in H$, $a\in A$.
\item[(ii)] $\sg$ is a cocycle. That is, $\sg(h,1)=\sg(1,h)=\vps(h)1$ for all $h\in H$ and
$$[h_1\cdot \sg(k_1,m_1)]\sg(h_2,k_2m_2)=\sg(h_1,k_1)\sg(h_2k_2,m)$$
for all $h,k,m\in H$.
\end{enumerate}
\end{lem}

Note that if $\sg$ is trivial, that is, $\sg(h,k)=\vps(h)\vps(k)1$, for all $h,k\in H$. Then the crossed product $A\#_\sg H$ is just the smash product $A\#H$.

\begin{rem}\label{cocykk}Let $A\#_\sg H$ be a crossed product and $\sg$ an invertible map in $\Hom(H\ot H, \kk)$. Then $A\#_\sg H$ is an associative algebra if and only if $\sg$ is a 2-cocycle and $A$ is an $H^\sg$-module algebra.
\end{rem}

\begin{defn}
Let $A\subseteq B$ be an extension of algebras, and $H$ a Hopf algebra.
\begin{enumerate}
\item[(i)] $A\subseteq B$ is called a \it{(right) $H$-extension} if $B$ is a right $H$-comodule algebra such that $B^{coH}=A$.
\item[(ii)] The $H$-extension $A\subseteq B$ is said to be \it{$H$-cleft} if there exists a right $H$-comodule morphism $\gm:H\ra B$ which is (convolution) invertible. Note that this $\gm$ can be chosen such that $\gm(1)=1$.
\end{enumerate}
If $\kk\subseteq B$ is $H$-cleft, then   $B$ is called a \it{(right) cleft object}. Left cleft extensions and left cleft objects can be defined similarly.
\end{defn}

\begin{lem}\label{mm}\cite[Theorem 7.2.2, Proposition 7.2.3, Proposition 7.2.7]{m} Let $H$ be a Hopf algebra. An $H$-extension $A\subseteq B$ is $H$-cleft with right convolution invertible $H$-comodule morphism $\gm:H\ra B$ if and only if $B\cong A\#_\sg H$ as algebras with a convolution invertible map $\sg:H\ot H\ra A$.   The
twisted $H$-module action on $A$ is given by
$$h\cdot a=\gm(h_1)a\gm^{-1}(h_2),$$
for all $a\in A$, $h\in H$.  Moreover, $\gm$ and $\sg$ are constructed each other by
$$\sg(h,k)=\gm(h_1)\gm(k_1)\gm^{-1}(h_2k_2)$$ and
$$\gm(h)=1\#h,\;\;\;\gm^{-1}(h)=\sg^{-1}(Sh_2,h_3)\#Sh_1$$
for all $h,k\in H$, $a\in A$.

\end{lem}

From this lemma, we see that right cleft objects of a Hopf algebra $H$ are just the algebras $_\sg H$, where $\sg$ is a 2-cocycle on $H$.

\subsection{AS-Gorenstein algebras}\label{s1as}

In this paper, unless otherwise stated, a graded algebra will always mean an $\NN$-graded algebra. An $\NN$-graded algebra $A =\op _{i\le 0}
 A_i$ is called connected if $A_0 = \kk$.

\begin{defn}

A connected graded algebra $A$ is called \it{AS-Gorenstein}
if the following conditions hold:
\begin{enumerate}
\item[(i)] $A$ has finite injective dimension $d$ on both sides,
\item[(ii)] $\Ext^i_A({}_A\kk, {}_AA) \cong \begin{cases}0,&i\neq d;\\ \kk(l),& i=d,\end{cases}$
where $l$ is an integer,
\item[(iii)] The right version of (ii) holds.
\end{enumerate}
If, in addtion,
\begin{enumerate}
\item[(iv)] $A$ is of  finite global dimension $d$,
then $A$ is called \it{AS-regular}.
\end{enumerate}
\end{defn}

Noe that an AS-Gorenstein (regular) algebra can be defined on an augmented algebra in general, see \cite{bz}.
For an algebra $A$, if the injective dimension of $_AA$ and $A_A$
are both finite, then these two integers are equal by \cite[Lemma
A]{za}. We call this common value the \textit{injective dimension} of $A$.
The left global dimension and the right global dimension of a
Noetherian algebra are equal.  When the global
dimension is finite, then it is equal to the injective dimension.

\begin{defn}(cf. \cite[defn. 1.2]{bz}). \label{defn as} Let $A$ be a Noetherian algebra with a fixed augmentation map $\vps:A\ra \kk$.
\begin{enumerate}
\item[(i)] The algebra $A$ is said to be \it{
AS-Gorenstein}, if
\begin{enumerate}
\item[(a)] $\injdim  {_AA}=d<\infty$,
\item[(b)] $\dim\Ext_A^i({_A\kk},{_AA})=\begin{cases}0,&i\neq d;\\ 1,& i=d,\end{cases}$
\item[(c)] the right versions of (a) and (b) hold.
\end{enumerate}

\item[(ii)] If, in addition, the global dimension of $A$ is finite, then $A$ is called  \it{AS-regular}.
\end{enumerate}
\end{defn}

The concept of a homological integral for an
AS-Gorenstein Hopf algebra was introduced by Lu,
Wu and Zhang in \cite{lwz} to study infinite dimensional Noetherian Hopf algebras. It is a generalization of the concept of an integral of a finite dimensional Hopf algebra. It turns out that homological integrals are useful in describing homological properties of Hopf algebras (see e.g. \cite[Theorem 2.3]{hvz}).

\begin{defn}\label{int}
 Let $A$ be an AS-Gorenstein algebra
with injective dimension $d$. Then $\Ext_A^d({_A\kk},{_AA})$ is a 1-dimensional
right $A$-module. Any nonzero element in $\Ext_A^d({_A\kk},{_AA})$
is called a \it{left homological integral} of $A$. We write
$\int_A^l$ for $\Ext_A^d({_A\kk},{_AA})$. Similarly, $\Ext_A^d({\kk_A},{A_A})$ is a 1-dimensional left $A$-module. Any nonzero element  in $\Ext_A^d({\kk_A},{A_A})$ is
called a \it{right homological integral} of $A$. Write
$\int_A^r$ for $\Ext_A^d({\kk_A},{A_A})$.

$\int_A^l$ and $\int_A^r$
are  called \textit{left and right homological integral modules} of $A$
respectively.
\end{defn}

The left integral module $\int^l_A$  is a 1-dimensional right
$A$-module. Thus $\int^l_A\cong \kk_\xi$ for some algebra
homomorphism $\xi:A\ra\kk$. Similarly,  $\int^r_A\cong {}_\eta\kk$ for some algebra homomorphism $\eta$.

\subsection{$N$-Koszul algebras}\label{s1k}

Let $V$ be a finite dimensional vector space, and $T(V)=\kk\ot V\ot V^{\ot 2}\ot \cdots$ be the tensor algebra with the usual grading.   A graded algebra $T(V)/\lan R\ran$ is called \it{N-homogenous} if $R$ is a subspace of $V^{\ot N}$.  Let $V^*$ be the dual space $\Hom(V,\kk)$. The algebra $A^!=T(V^*)/\lan R^{\perp}\ran$ is called the \it{homogeneous dual} of $A$, where $R^{\perp}$ is the orthogonal subspace of $R$ in $(V^*)^{\ot N}$.

\begin{rem}
Let $\phi$ be the map defined as follows:
$$\begin{array}{rccl}
\phi:&(V^*)^{\ot n}&\ra& (V^{\ot n})^*\\
&f_n\ot f_{n-1}\ot \cdots \ot f_1&\mapsto& \phi(f_n\ot f_{n-1}\ot \cdots \ot f_1),
\end{array}$$
where $\phi(f_n\ot f_{n-1}\ot \cdots \ot f_1)(x_1\ot \cdots x_{n-1}\ot x_n)=f_1(x_1)f_2(x_2)\cdots f_n(x_n)$, for any $x_1\ot x_2\ot \cdots\ot x_n\in V^{\ot n}$. This map $\phi$ is a bijection.
Throughout, we identify $(V^*)^{\ot n}$ with $(V^{\ot n})^*$ via this bijection.
\end{rem}

Let ${\bf n}:\NN\ra \NN$ be the function defined by
$${\bf n}(i)=\begin{cases}Nk,&i=2k\\Nk+1,&i=2k+1.\end{cases}$$

An $N$-homogenous algebra $A$ is called \it{$N$-Koszul} if the trivial module ${}_A\kk$ admits a graded projective resolution
$$\cdots \ra P_i\ra P_{i-1}\ra \cdots \ra P_1\ra P_0\ra {}_A\kk\ra 0$$
such that $P_i$ is generated in degree ${\bf n}(i)$ for all $i\le 0$. A Koszul algebra is a just 2-Koszul algebra.

The Koszul bimodule complex of a Koszul algebra is constructed by Van den Bergh in \cite{vdb2}. This complex was generalized to  $N$-Koszul case in \cite{bm}.
Now let $A=T(V)/\lan R\ran$ be an $N$-Koszul algebra. Let $\{e_i\}_{i=1,2,\cdots,n}$ be a basis of $V$ and $\{e^*_i\}_{i=1,2,\cdots,n}$ the dual basis.  Define two $N$-differentials $$d_l,d_r:A\ot (A^!_p)^*\ot A\ra A\ot (A^!_{p-1})^*\ot A$$ as follows:
\begin{eqnarray*}
& d_l(x\ot \om \ot y)=\sum_{i=1}^n x e_i\ot  e_i^*\cdot\om \ot y\\
& d_r(x\ot \om \ot y)=\sum_{i=1}^nx \ot \om \cdot e_i^*\ot e_iy,
\end{eqnarray*}
for $x\ot \om \ot y\in A\ot (A^!_p)^*\ot A$. The left action $e_i^*\cdot\om$ is defined by $[e_i^*\cdot\om](\al)=\om(\al e_i^*)$ for any $\al\in (A^!_{p-1})^*$. The right action $\om \cdot e_i^*$ is defined similarly. One can check that $d_l$ and $d_r$ commute. Fix a primitive $N$-th root of unity $q$. Define $d:A\ot (A^!_p)^*\ot A\ra A\ot (A^!_{p-1})^*\ot A$ by $d=d_l-q^{p-1}d_r$. We obtain the following $N$-complex:
{\small$$\mathbf{K_{l-r}(A)}:\cdots \xra{d_l-d_r} A\ot (A^!_N)^*\ot A\xra{d_l-q^{N-1}d_r} \cdots \xra{d_l-qd_r} A\ot V\ot A\xra{d_l-d_r} A\ot A\ra 0.$$}
The bimodule Koszul complex $\mathbf{K_{b}(A)}$ is a contraction of $\mathbf{K_{l-r}(A)}$. It is obtained by keeping the arrow $A\ot V\ot A\xra{d_l-d_r}A\ot A$ at the far right, then putting together the $N-1$ consecutive ones, and continuing alternately:
{\small$$\mathbf{K_{b}(A)}:\cdots \xra{d^{N-1}} A\ot (A^!_{N+1})^*\ot A\xra{d}A\ot (A^!_{N})^*\ot A \xra{d^{N-1}} A\ot V\ot A\xra{d} A\ot A\ra 0.$$}
Here $d=d_l-d_r$ and $d^{N-1}=d_l^{N-1}+d_l^{N-2}d_r+\cdots+d_ld_r^{N-2}+d_r^{N-1}$.

An $N$-homogenous algebra is $N$-Koszul if and only if the complex $\mathbf{K_{b}(A)}\ra A\ra 0$ is exact via the multiplication $A\ot A\ra A$ \cite[Theorem 4.4]{bm}. Moreover, in such a case, $\mathbf{K_{b}(A)}\ra A\ra 0$ is a minimal bimodule free resolution of $A$.

\subsection{Calabi-Yau algebras}\label{s1ca}
\begin{defn}\label{tcy}An algebra $A$ is called  a \it{twisted Calabi-Yau algebra  of dimension
$d$} if
\begin{enumerate}\item[(i)] $A$ is \it{homologically smooth}, that is, $A$ has
a bounded resolution of finitely generated projective
$A^e$-modules; \item[(ii)] There is an automorphism $\mu$ of $A$ such that
\begin{equation}\label{cy1}\Ext_{A^e}^i(A,A^e)\cong\begin{cases}0,& i\neq d
\\A^\mu,&i=d\end{cases}\end{equation}
as $A^e$-modules.
\end{enumerate}
If such an automorphism $\mu$ exists, it is unique up to an inner automorphism and is called the \it{Nakayama automorphism} of $A$. A \it{Calabi-Yau algebra} is a twisted Calabi-Yau algebra whose Nakayama automorphism is an inner automorphism.
\end{defn}

A \it{Graded twisted CY algebra} can be defined in a similar way. That is, we should consider the category of graded modules and condition (\ref{cy1}) should be replaced by
$$\Ext_{A^e}^i(A,A^e)\cong\begin{cases}0,& i\neq d;
\\A_\mu(l),&i=d,\end{cases}$$
where $l$ is an integer and $A_\mu(l)$ is the shift of $A_\mu$ by degree $l$.

We end this section with the following lemma, which shows that AS-regular Hopf algebras are just twisted CY Hopf algebras.
\begin{lem}\label{cyas}
Let $A$ be a Noetherian AS-regular Hopf algebra with $\int_A^l=\kk_\xi$, where $\xi:A\ra\kk$ is an algebra homomorphism. The followings hold:
\begin{enumerate}
\item[(i)]\cite[Lemma 1.3]{rrz} The algebra $A$ is twisted CY with Nakayama automorphism $\mu$ defined by $\mu(x)=S^{-2}(x_1)\xi(x_2)$ for any $x\in A$. (Alternatively, the algebra automorphism $\nu$ defined by $\nu(x)=\xi(x_1)S^2(x_2)$ is
also a Nakayama automorphism of $A$).
\item[(ii)]\cite[Theorem 2.3]{hvz} The algebra $A$ is CY if and only if $\xi=\vps$, and $S^2$ is an inner automorphism.
\end{enumerate}
\end{lem}

\section{The CY property of Cleft extension}\label{s2}

Let $H$ be a Hopf algebra, $\sg$ a 2-cocycle on $H$ and $A$ an $N$-Koszul $H^\sg$-module algebra. Then the crossed product $A\#_\sg H$ is an associative algebra. In this section we  show that $A\#_\sg H$ is a twisted CY algebra if both $A$ and $H$ are twisted CY algebras. This generalizes \cite[
Theorem 2.12]{liwz} and \cite[Theorem 0.2]{rrz},

The following definition is inspired by ``$H_{S^i}$-equivariant $A$-bimodule'' introduced in \cite[Definition 2.2]{rrz}, where $H$ is a Hopf algebra and $i$ is an even integer.
\begin{defn}\label{def}
Let $H$ be a Hopf algebra and $A$ a left $H$-module algebra. For a given even integer $i$, we define an algebra $A^e\rtimes_{S^i}H$. As vector spaces, $A^e\rtimes_{S^i}H=A\ot A\ot H$. The multiplication is given by
$$(a\ot b\ot g)(a'\ot b'\ot h)=a(S^ig_1\cdot a')\ot (g_3\cdot b')b\ot g_2h,$$
\end{defn}
for any $a\ot b\ot g, a'\ot b'\ot h\in A\ot A\ot H$.

\begin{rem}
\begin{enumerate}
\item[(i)] When $i=0$, $A^e\rtimes_{S^i}H$ is just the algebra $A^e\rtimes H$ introduced  by Kaygun \cite{kay}.

\item[(ii)]  An $A^e\rtimes_{S^i}H$-module $M$ is a vector space such that it is both an $A^e$-module and an $H$-module satisfying
\begin{equation}\label{ash}h\cdot (amb)=((S^ih_1)\cdot a) (h_2\cdot m)(h_3\cdot b),\end{equation}for any $h\in H$, $a, b\in A$ and $m\in M$.
\end{enumerate}

\end{rem}

\begin{lem}\label{hhom1}Let $M$ be an $A^e\rtimes_{S^i} H$-module and $N$ an $(A\#H)^e$-module.
\begin{enumerate}
\item[(i)] The space $\Hom_{A^e}(M,N)$ is a left $H$-module with the $H$-action defined by
\begin{equation}\label{hhom}(h\rightharpoonup f)(m)=(S^ih_3)f[(S^{-1}h_2)\cdot m](S^{-1}h_1)\end{equation}
for any $h\in H$, $f\in \Hom_{A^e}(M,N)$ and $m\in M$.
\item[(ii)] The space $M\ot_{A^e}N$ is a left $H$-module with the $H$-action given by
\begin{equation}\label{hten}h\cdot (m\ot n)=h_2\cdot m\ot h_3n (S^{i+1}h_1)\end{equation}
for any $h\in H$ and $m\ot n\in M\ot  N$.
\end{enumerate}
\end{lem}
\proof The proof is routine and quite similar to the proofs of Lemma 1.8 and Lemma 1.9 in \cite{liwz}.

\begin{rem}
Keep the notations as in Lemma \ref{hhom1}, $\Hom_{A^e}(M,N)$ can be made into a right $H$-module by defining $f\leftharpoonup h=Sh\rightharpoonup f$ for any $h\in H$ and $f\in \Hom_{A^e}(M,N)$. That is,
\begin{equation}\label{rha}(f\lap h)(m)=S^{i+1}h_1f(h_2\cdot m)h_3.\end{equation}
\end{rem}

Since $A$ is a left $H$-module algebra, the algebra $A^e$ is an $(A\#H)^e$-module with the following module structure:
\begin{equation}\label{ae}(a\#h)\cdot(x\ot y)=a(h\cdot x)\ot y, \;\;
(x\ot y)\cdot(b\# g)=x\ot (S^{-1}g)\cdot (yb)\end{equation}
for any $x\ot y\in A^e$ and $a\#h$, $b\# g\in A\#H$.

By Lemma \ref{hhom1}, $\Hom_{A^e}(M,A^e)$ is a left $H$-module for any $A^e\rtimes H$-module $M$. Furthermore, the $A^e$-bimodule structure of $A^e$ induces a left $A^e$-module structure on $\Hom_{A^e}(M,A^e)$. That is,
\begin{equation}\label{maa}[(a\ot b)\cdot f](x)=f(x)(b\ot a),\end{equation}
for any $a\ot b\in A^e$, $f\in \Hom_{A^e}(M,A^e)$ and $x\in M$.

In \cite{liwz} the authors showed that if $H$ is involutory, then $\Hom_{A^e}(M,A^e)$ is again an $A^e\rtimes H$-module for any $A^e\rtimes H$-module $M$. In general, we have the following.

\begin{lem}\label{s-2}
Let $M$ be an $A^e\rtimes H$-module. Then $\Hom_{A^e}(M,A^e)$ is an \\ $A^e\rtimes_{S^{-2}}H$-module.
\end{lem}

%\begin{rem}\label{rem}$\Hom_{A^e}(M,A^e)$ is a right $H$-module as defined in (\ref{rha}). This lemma shows that the right $H$-module structure and the $A^e$-module structure satisfies
%\begin{equation}\label{comp} (afb)\leftharpoonup h=(S^{-1}h_3\cdot a)(f \leftharpoonup h_2)(Sh_1\cdot b)\end{equation}
%for all $h\in H$, $a,b\in A$ and $f\in \Hom_{A^e}(M,A^e)$.
%\end{rem}

In \cite[Theorem 2.4]{liwz} the Van den Bergh duality was generalized  to algebras with a Hopf action from an involutory Hopf algebra. In fact, we can drop the condition ``involutory''.

\begin{prop}
Let $H$ be a Hopf algebra and $A$ a left  $H$-module algebra. Assume that $A$ admits a finitely generated $A^e$-projective resolution of
finite length such that it is a complex of  $A^e\rtimes H$-modules. Suppose there exists an integer $d$ such that
$$\Ext^i_{A^e} (A, A^e)=\begin{cases}0, &i\neq d;\\
U,&i=d,
\end{cases}$$
where $U$ is an invertible $A^e$-module.
Then for any $(A\#H)^e$-module
$N$, we have
$$\HH^i(A, N ) \cong {}_{S^{-2}}\HH_{d-i} (A, U\ot_A N)$$
as left $H$-modules.
\end{prop}
\proof Suppose that $P$ is an $A^e\rtimes H$-module such that it is finitely generated and projective as an $A^e$-module, and $N$ is an $(A\#H)^e$-module.
By Lemma \ref{s-2}, $\Hom_{A^e}(P,A^e)$ is an $A^e\rtimes_{S^{-2}} H$-module.   So $\Hom_{A^e}(P,A^e)\ot_{A^e}N$ is an $H$-module with the module structure given by (\ref{hten}). Moreover, the equation (\ref{hhom}) defines an $H$-module structure on $\Hom_{A^e}(P,N)$. With these $H$-actions, one can check that the canonical isomorphism
 $$\Psi:\Hom_{A^e}(P,A^e)\ot_{A^e}N\ra \Hom_{A^e}(P,N)$$ is also an $H$-isomorphism. Therefore, the proof of \cite[Theorem 2.4]{liwz} works for non-involutory Hopf algebras. But for a non-involutory Hopf algebra $H$,  the module $U$ is an $A^e\rtimes_{S^{-2}} H$-module by Lemma \ref{s-2}. Thus,
$U\ot_A N$ is an $(A\#H)^e$-module with module structure defined by
\begin{equation}\label{un}(a\#h)\cdot (u\ot n)=a((S^2h_1)\cdot u)\ot (S^2h_2)\cdot n,\;\;(u\ot n)\cdot(b\#g)=u\ot n\cdot (b\#g),\end{equation}
for any $a\#h$, $b\#g\in A\#H$ and $u\#n\in U\ot N$. Consequently, we have the following $H$-isomorphisms:
$$\begin{array}{ccl}\HH^i(A,N)&\cong&\Ext_{A^e}^i(A,N)\\
&\cong& \H^i(\RHom_{A^e}(A,N))\\
&\cong& \H^i(\RHom_{A^e}(A,A^e){}^L\ot_{A^e}N)\\
&\cong& \H^i(U[-d]{}^L\ot_{A^e}N)\\
&\cong& \H^{i-d}(U{}^L\ot_{A^e}N)\\
&\cong& \H^{i-d}({}_{S^{-2}}[A\ot_{A^e}(U{}^L\ot_AN)])\\
&\cong& {}_{S^{-2}}\HH_{d-i}(A,U\ot_AN).\\
\end{array}
$$\qed

In the rest of this section, we work with the category of graded modules. Let $A$ be a graded algebra, and let  $A$-GrMod denote the category of graded left $A$-modules and graded homomorphisms of degree zero. For any $M,N\in A$-GrMod, $\Hom_A(M,N)$ is the graded vector space consisting of graded $A$-module homomorphisms. That is,
$$\Hom_A(M,N)=\op_{i\in \ZZ}\Hom_{A\t{-GrMod}}(M,N(i)).$$

Let $H$ be a Hopf algebra. We say that a graded algebra $A$ is a left graded $H$-module algebra if it is a left $H$-module algebra such that each $A_i$ is an $H$-module. Let $\sg$ is a 2-cocycle on $H$. The cocycle deformation $H^\sg$ is a Hopf algebra. If $A$ is a left graded $H^\sg$-module algebra, then we have the algebra $A\#H^\sg$. Moreover, we can construct the algebra $A\#_\sg H$ by Remark \ref{cocykk}.  It is easy to see that both $A\#H^\sg$ and $A\#_\sg H$ have natural graded algebra structures.

Now, we fix a Hopf algebra $H$ and a 2-cocycle $\sg$ on $H$. Let $V$ be a left $H^\sg$-module and $A=T(V)/\lan R\ran$ an $N$-Koszul graded $H^\sg$-module algebra. The dual $V^*$ is a right $H^\sg$-module with the module structure given by
\begin{equation}\label{lhd}
(\al\lhd h)(x)=\al (h\cdot x).\end{equation}
for  $\al \in V^*$, $h\in H$ and $x\in V$.
\begin{rem}\label{**}
Let $\{e_1,e_2,\cdots,e_n\}$ be a basis of $V$. Suppose that $h\cdot e_i=\sum_{j=1}^nc^h_{ji}e_j$ with $c^h_{ji}\in \kk$. Then we have  $e_i^*\lhd h=\sum_{j=1}^nc^h_{ij}e_j^*$.
\end{rem}

We extend the action ``$\lhd$'' on $V^*$ to $(V^*)^{\ot n}$:
$$(\al_n\ot \al_{n-1}\ot \cdots \ot \al_1)\lhd h=(\al_n \lhd h_n)\ot (\al_{n-1} \lhd h_{n-1})\ot \cdots \ot (\al_1\lhd h_1).$$

It is easy to check that $R^\perp \lhd h\subseteq R^\perp$. Consequently, $A^!$ is a right $H^\sg$-module algebra with the action ``$\lhd$''. In fact, one can make $A^!$ into a left $H^\sg$-module algebra as follows:
\begin{equation}\label{ha1}h\cdot \bt =\bt \lhd (S^{\sg^{-1}} h), \end{equation}
for any $\bt\in A^!$ and $h\in H$.

Thanks to Lemma \ref{s-2}, we obtain the following proposition generalizing  \cite[Proposition 2.2]{liwz}.

\begin{prop}\label{ext}
Let $H$ be a Hopf algebra, $\sg$ a 2-cocycle on $H$, and $A$ a left graded $H^\sg$-module algebra. If $A$ is an $N$-Koszul graded twisted CY algebra of dimension $d$ with Nakayama automorhism $\mu$, then as $A^e\rtimes_{S^{ \sg^{-2}}} H^\sg$-modules
$$\Ext_{A^e}^i(A,A^e)\cong \begin{cases}0, &i\neq d;\\A_{\mu}\ot A^!_{{\bf n}(d)}, &i=d,\end{cases}$$
where the $A^e\rtimes _{S^{\sg^{-2}}} H^\sg$-module  structure on $A_\mu\ot A^!_{{\bf n}(d)}$ is given by
\begin{equation}\label{a1d}(a\ot b\ot h)(x\ot \al)=a((S^{\sg^{-2}}h_1)\cdot x)\mu(b)\ot h_2\cdot \al,\end{equation}
for any $a\ot b\ot h\in A^e\rtimes_{S^{-2}}  H$ and $x\ot \al\in A\ot A^!_{{\bf n}(d)}$.
\end{prop}
\proof The algebra $H^\sg$ is a Hopf algebra and the algebra $A$ is a left $H^\sg$-module algebra. Proposition 2.1 in \cite{liwz} shows that the $A^e$-projective resolution $\mathbf{K_b(A)}\ra A\ra0$ of $A$ is
an $A^e\rtimes H^\sg$-module complex. The $A^e\rtimes H^\sg$-module structure is defined as follows. Each term in $\mathbf{K_b(A)}$ is of the form $A\ot (A_p^!)^*\ot A$.  Since  $A^!_p$ is a right $H^\sg$-module with the action ``$\lhd$'' defined in (\ref{lhd}), $(A_p^!)^*$ is a natural left $H^\sg$-module. That is,
$$(h\cdot \om)(x)=\om(x\lhd h),$$
for any $h\in H^\sg$, $\om\in (A_p^!)^*$ and $x\in A_p^!$. Each $A\ot (A_p^!)^*\ot A$ is an $A^e\rtimes H^\sg$-module with the module structure defined by
\begin{equation}\label{hk}(a\ot b\ot h)\cdot(x\ot \om\ot y)=a( h_1\cdot x)\ot h_2\cdot \om \ot (h_3 \cdot y)b,\end{equation}
where $a\ot b\ot h \in A^e\rtimes H$ and $ x\ot \om\ot y \in A\ot (A_p^!)^*\ot A$.

Now we recall another bimodule complex constructed in \cite{bm}. First, we define two $N$-differentials:
$$\dt_l,\dt_r:A\ot A^!_p\ot A\ra A\ot A^!_{p+1}\ot A$$ as follows:
$$\dt_l(x\ot \al \ot y)= \sum_{i=1}^nx e_i\ot  e_i^*\al \ot y, \mathrm{and}\
\dt_r(x\ot \al \ot y)=\sum_{i=1}^nx \ot \al e_i^*\ot e_iy,$$
 for $x\ot \al \ot y\in A\ot A^!_p\ot A$. It is easy to check that $\dt_l$ and $\dt_r$ commute. Fix a primitive $N$-th root of unity $q$. The complex
{\small $$ \mathbf{L_{l-r}(A)}:A\ot A\xra{\dt_r-\dt_l}  A\ot V^*\ot A\xra{\dt_r-q\dt_l}\cdots\xra{\dt_r-q^{N-1}\dt_l} A\ot A^!_{N}\ot A\xra{\dt_r-\dt_l}\cdots$$}
is an $N$-complex. The complex $\mathbf{L_b(A)}$ is the contraction of $ \mathbf{L_{l-r}(A)}$. It is obtained by keeping the arrow $A\ot A\xra{\dt_r-\dt_l}  A\ot V^*\ot A$ at the far left, then putting together the $N-1$ following ones, and continuing alternately:
{\small $$\mathbf{L_b(A)}:A\ot A\xra{\dt}  A\ot V^*\ot A\xra{\dt^{N-1}} A\ot A^!_N\ot A\xra{\dt} A\ot A^!_{N+1}\ot A\xra{\dt^{N-1}} \cdots,$$}
where $\dt=\dt_r-\dt_l$ and $\dt^{N-1}=\dt_r^{N-1}+\dt^{N-2}_r\dt_l+\cdots+\dt_r\dt_l^{N-2}+\dt_l^{N-1}$.
When the Hopf algebra  $H^\sg$ is involutory, Proposition 2.2 in \cite{liwz} shows that the complex $\Hom_{A^e}(\mathbf{K_b(A)},A^e)$ and the complex $\mathbf{L_b(A)}$ are isomorphic as $A^e\rtimes H^\sg$-complexes.

When $H^\sg$ is not involutory, $\Hom_{A^e}(\mathbf{K_b(A)},A^e)$ is a complex of $A^e\rtimes_{S^{\sg^{-2}}} H^\sg$-modules by Lemma \ref{s-2}. In this case,
 $\Hom_{A^e}(\mathbf{K_b(A)},A^e)$ and $\mathbf{L_b(A)}$ are isomorphic as $A^e\rtimes_{S^{\sg^{-2}}} H^\sg$-module complexes.  The $A^e\rtimes_{S^{\sg^{-2}}} H^\sg$-module structure of each term $A\ot A^!_p\ot A$ in $\mathbf{L_b(A)}$ is  given by
$$(a\ot b\ot h )\cdot (x\ot \al \ot y)=a((S^{\sg^{-2}}h_1)\cdot x)\ot h_2\cdot \al \ot (h_3\cdot y)b,$$
for any $a\ot b\ot h\in A^e\rtimes_{S^{\sg^{-2}}} H^\sg$ and $x\ot \al \ot y\in A\ot A^!_p\ot A$.

Now we can use the complex $\mathbf{L_b(A)}$ to compute $\Ext^*_{A^e}(A,A^e)$. The method is the same as the one in the proof of Proposition 2.2 in \cite{liwz}.

Since the algebra $A$ is an $N$-Koszul graded twisted CY algebra,  $A$ is AS-regular (see \cite[Lemma 1.2]{rrz}). The Ext algebra $E(A)$ of $A$ is graded Frobenius by Corollary 5.12 in \cite{bm}. Thus, there exists an automorphism $\phi$ of $E(A)$, such that
$$E(A)_\phi\cong E(A)^{*}(-d)$$
as $E(A)$-bimodules.

Let $\{e_1,e_2,\cdots,e_n\}$ be a basis of $A_1=V$, and $\{e_1^*,e_2^*,\cdots,e_n^*\}$ the corresponding dual basis.
Suppose that $\phi$ is given by $$\phi(e_1^*,e_2^*,\cdots,e_n^*)=(e_1^*,e_2^*,\cdots,e_n^*)Q,$$
for some invertible matrix $Q$. Define an automorphism of $A$ via $$\vph(e_1,e_2,\cdots,e_n)=(e_1,e_2,\cdots,e_n)Q^T,$$
where $Q^T$ is the transpose of $Q$. It is obvious that the restriction of $\phi$ to $V^*$ and the restriction of $\vph$ to $V$ are dual to each other.

Let $\epsilon$ be the automorphism of $A$ defined by  $\epsilon(a)=(-1)^ia$ for any homogeneous element $a\in A_i$.  By assumption, we have $\Ext_{A^e}^i(A,A^e)=0$ for $i\neq d$. Now we compute $\Ext^d_{A^e}(A, A^e)$. Suppose $N\le 3$. Then the  dimension $d$ must be odd. We consider the following sequence
\begin{equation}\label{com}A\ot A_{{\bf n}(d)-1}^!\ot A\xra{\dt} A\ot A_{{\bf n}(d)}^!\ot A\xra{u} A_{\mu}\ot A_{\mathbf{ n}(d)}^!\ra 0,\end{equation}
where $\mu=\eps^{d+1}\vph$ and the morphism $u$ is given by $u(x\ot \al\ot y)=x\mu(y)\ot \al$, for any $x\ot \al\ot y\in A\ot A_{\mathbf{ n}(d)}^!\ot A$. Since $E(A)$ is Frobenius with Nakayama automorphism $\phi$, by \cite[Proposition 3.1]{bm},  we have $e_i^*\al=\al\phi(e_i^*)$, for any $\al \in
A^!_{{\bf n}(d)-1}$. Now for any $x\ot \al\ot y\in A\ot A^!_{\mathbf{ n}(d)-1}\ot A$, we have:
{\small $$\begin{array}{ccl}u\dt(x\ot \al\ot y)
&=&u(\sum_{i=1}^nx\ot \al e_i^*\ot e_iy-\sum_{i=1}^nxe_i\ot e_i^*\al\ot y)\\
&=&\sum_{i=1}^nx\mu(e_iy)\ot \al e_i^*-\sum_{i=1}^nxe_i\mu(y)\ot e_i^*\al\\
&=&\sum_{i=1}^nx\mu(e_i)\mu(y)\ot \al e_i^*-\sum_{i=1}^nxe_i\mu(y)\ot \al \phi(e_i^*)\\
&=&\sum_{i=1}^nx\mu(e_i)\mu(y)\ot \al e_i^* -\sum_{i=1}^nxe_i\mu(y)\ot \al (\sum_{j=1}^nq_{ji}e^*_j) +\\
&=&\sum_{i=1}^nx\mu(e_i)\mu(y)\ot \al e_i^*-\sum_{i=1}^n\sum_{j=1}^nq_{ji}xe_i\mu(y)\ot \al e^*_j \\
&=&\sum_{i=1}^nx\mu(e_i)\mu(y)\ot \al e_i^*-\sum_{i=1}^n x\vph(e_i)\mu(y)\ot \al e^*_i \\
&=&\sum_{i=1}^n(-1)^{d+1}x\vph(e_i)\mu(y)\ot \al e_i^*-\sum_{i=1}^n x\vph(e_i)\mu(y)\ot \al e^*_i \\
&=&0.
\end{array}$$}
Therefore, the sequence (\ref{com}) is  a complex. Hence, it is exact by \cite[Proposition 4.1]{bm}.

Similar to the proof of \cite[Prop 2.2]{liwz}, we can show that
 (\ref{a1d}) defines an $A^e\rtimes_{S^{\sg^{-2}}} H^\sg$-module structure on $A\ot A_d^!$ and $u$ is an $A^e\rtimes_{S^{\sg^{-2}}} H^\sg$-homomorphism. Therefore, $\Ext_{A^e}^d(A,A^e)\cong A_\mu\ot A^!_{{\bf n}(d)}$ as $A^e\rtimes_{S^{\sg^{-2}}} H^\sg$-modules.

 For the case $N=2$, the proof is similar.
\qed

Let $H$ be a Hopf algebra, $\sg$ a 2-cocycle on $H$, and $A$ a graded $H^\sg$-module algebra. Let $P$ be an $A^e\rtimes H^\sg$-module. $\Hom_{A^e}(P,A^e)$ is a right $H^\sg$-module as defined in (\ref{rha}). Then we can define a right $H$-module structure on $\Hom_{A^e}(P,A^e)\ot {}_\sg H\ot {}_\sg H$:
\begin{equation}\label{h1}(f\ot k\ot l)\leftarrow h=f\leftharpoonup h_2\ot (S_{1,\sg}h_1)\centerdot_\sg k\ot l \centerdot_\sg h_3\end{equation}
for all $f\ot k\ot l\in \Hom_{A^e}(P,A^e)\ot {}_\sg H\ot {}_\sg H$ and $h\in H$. Recall that $H$ can be viewed as the algebra $H(1,1)$. Here $h_1\ot h_2\ot h_3=(\bdt_{1,\sg}^\sg\ot \id)\bdt_{1,1}^\sg(h)$.  Both $\bdt_{1,\sg}^\sg$ and $\bdt_{1,1}^\sg$ are algebra homomorphisms. So this $H$-module is well-defined.  We denote this $H$-module by $\Hom_{A^e}(P,A^e)_*\ot {}_{*\sg} H\ot {}_\sg H_*$.

The right $H$-module structure of $H$ induces a natural $H$-module structure on  $\Hom_{A^e}(P,A^e)\ot {}_\sg H\ot H$. That is,
\begin{equation}\label{h2}(f\ot k\ot l)\leftarrow h =f\ot k\ot lh\end{equation}
for all $f\ot k\ot l\in \Hom_{A^e}(P,A^e)\ot {}_\sg H\ot  H$ and $h\in H$.  We denote this $H$-module by $\Hom_{A^e}(P,A^e)\ot {}_{\sg} H\ot  H_*$.

We can define an $(A\#_\sg H)^e$-module structure on $\Hom_{A^e}(P,A^e)\ot {}_\sg H\ot H _*$ as follows:
\begin{equation}\label{phh1}\begin{array}{l}(a\#h)\cdot(f\ot k\ot l)=a((S^{\sg^2} h_1)\rightharpoonup f )\ot S_{1,\sg}(S_{\sg,1}(h_2))\ced_\sg k\ot h_3 l,\\
(f\ot k\ot l)\cdot(b\#g)=f(k_1\cdot b)\ot k_2\ced_\sg g\ot l,\end{array}\end{equation}
for any $a\#h$, $b\#g\in A\#_\sg H$ and $f\ot k\ot l\in \Hom_{A^e}(P,A^e)\ot {}_\sg H\ot H_*$. Recall that the left $H^\sg$-module structure of $\Hom_{A^e}(P,A^e)$ is defined in (\ref{hhom}). Here $h_1\ot h_2\ot h_3=(\bdt_{\sg,1}^\sg\ot \id)\bdt_{\sg,1}^\sg(h)$ and $k_1\ot k_2=\bdt_{\sg,1}^\sg(k)$. We first check that the left $A\#_\sg H$-module structure is well-defined. We have the following equations:
$$\small{\begin{array}{rl}
&(b\#g)\cdot[(a\#h)\cdot(f\ot k\ot l)]\\
=&(b\#g)\cdot[a((S^{\sg^2} h_1)\rightharpoonup f)\ot S_{1,\sg}(S_{\sg,1}(h_2))\ced_\sg k\ot h_3 l]\\
=&b[(S^{\sg^2} g_1)\rightharpoonup(a((S^{\sg^2} h_1)\rightharpoonup f))]\ot S_{1,\sg}(S_{\sg,1}(g_2))\ced_\sg S_{1,\sg}(S_{\sg,1}(h_2))\ced_\sg k\ot g_3 h_3 l\\
\overset{(\ref{ash})}=&b[(g_1\cdot a)((S^{\sg^2} g_2)*(S^{\sg^2} h_1))\rightharpoonup f]\ot S_{1,\sg}(S_{\sg,1}(g_3))\ced_\sg S_{1,\sg}(S_{\sg,1}(h_2))\ced_\sg k\ot g_4 h_3 l\\
=&b(g_1\cdot a)((S^{\sg^2} g_2)*(S^{\sg^2} h_1))\rightharpoonup f\ot S_{1,\sg}(S_{\sg,1}(g_3))\ced_\sg S_{1,\sg}(S_{\sg,1}(h_2))\ced_\sg k\ot g_4 h_3 l\\
=&b(g_1\cdot a)(S^{\sg^2} (g_2* h_1))\rightharpoonup f\ot S_{1,\sg}(S_{\sg,1}(g_3\ced_\sg h_2))\ced_\sg k\ot g_4 h_3 l\\
=&[b(g_1\cdot a)\# g_2\ced_\sg h]\cdot(f\ot k\ot l)\\
=&[(b\#g)(a\#h)]\cdot(f\ot k\ot l).
\end{array}}$$
By Lemma \ref{anpo} we know that $S_{1,\sg}\circ S_{\sg,1}$ is an algebra homomorphism of ${}_\sg H$. Therefore, the fifth equation holds. The sixth equation follows from the fact that $\bdt_{\sg,1}^\sg$ is an algebra homomorphism. It follows that $\Hom_{A^e}(P,A^e)\ot {}_\sg H\ot H _*$ is a left $A\#_\sg H$-module. Similarly, we can see that $\Hom_{A^e}(P,A^e)\ot {}_\sg H\ot H _*$ is a right $A\#_\sg H$-module and for any $a\#h, b\#g\in A\#_\sg H$, and $f\ot k\ot l\in \Hom_{A^e}(P,A^e)\ot {}_\sg H\ot H _*$,
$$[(a\#h)(f\ot k\ot l)](b\#g)=(a\#h)[(f\ot k\ot l)(b\#g)].$$ In conclusion, $\Hom_{A^e}(P,A^e)\ot {}_\sg H\ot H _*$ is indeed an $(A\#_\sg H)^e$-module as defined in (\ref{phh1}).

The module $\Hom_{A^e}(P,A^e)_*\ot {}_{*\sg}H\ot {}_\sg H_*$ is also an $(A\#_\sg H)^e$-module with the module structure defined by
\begin{equation}\label{phh2}\begin{array}{l}(a\#h)\cdot(f\ot k\ot l)=(S^{\sg^{-1}}(h_1l_1)\cdot a)f\ot k\ot h_2\ced_\sg l_2,\\
(f\ot k\ot l)\cdot(b\#g)=f(k_1\cdot b)\ot k_2\ced_\sg g\ot l,\end{array}\end{equation} where $h_1\ot h_2=\bdt_{\sg,1}^\sg(h)$, $l_1\ot l_2=\bdt_{\sg,1}^\sg(l)$ and $k_1\ot k_2=\bdt_{\sg,1}^\sg(k)$.

Now both $\Hom_{A^e}(P,A^e)_*\ot {}_{*\sg}H \ot H_*$ and $\Hom_{A^e}(P,A^e)\ot  {}_\sg H\ot H _*$ are right $H\ot (A\#_\sg H)^e$-modules.

\begin{lem}\label{vh1}
Let $H$ be a Hopf algebra, $\sg$ a 2-cocycle on $H$, and $A$ a graded left $H^\sg$-module
algebra. If $P$ is an $A^e\rtimes H^\sg$-module, then the following $\Psi$ and $\Phi$ are $H\ot (A\#_\sg H)^e$-module isomorphisms
$$\Hom_{A^e}(P,A^e)_*\ot {} _{*\sg}H\ot {}_\sg H_*\overset{\Psi}{\underset{\Phi}\rightleftarrows}\Hom_{A^e}(P,A^e)\ot H\ot H_*,$$
where the module structures are given by (\ref{h1}), (\ref{h2}), (\ref{phh1}) and (\ref{phh2}),  $\Psi$ and $\Phi$ are defined as follows:
$$\begin{array}{ccl}\Psi(f\ot k\ot l)&=& f \leftharpoonup S^\sg(l_1)\ot S_{1,\sg}(S_{\sg,1}(l_2))\ced_\sg k\ot l_3,\\
\Phi(f\ot k\ot l)&=& f \leftharpoonup l_2\ot S_{1,\sg}(l_1)\ced_\sg k\ot l_3.\end{array}$$
Moreover, $\Psi$ and $\Phi$ are inverse to each other.
\end{lem}

\begin{lem}\label{rh}
Let $H$ be a Hopf algebra, $\sg$ a 2-cocycle on $H$, and $A$ a graded left $H^\sg$-module
algebra.  Let $P$ be an $A^e\rtimes H^\sg$-module, and $M$ an $(A\#_\sg H)^e$-bimodule. Then $\Hom_{A^e}(P,M)$ is a right $H$-module defined by
$$(f\leftarrow h)(x)=S_{1,\sg}(h_1)f(h_2x)h_3$$
for any $h\in H$, $f\in \Hom_{A^e}(P,M)$ and $x\in P$. Here $h_1\ot h_2\ot h_3=(\Delta_{1,\sg}^\sg\ot \id)\bdt_{1,1}^\sg(h)$.
\end{lem}
\proof For any $h,k\in H$ and $f\in \Hom_{A^e}(P,M)$, the following equations hold:
$$\begin{array}{ccl}[(f\leftarrow h)\la k](x)&=&S_{1,\sg}(k_1)(f\leftarrow h)(k_2x)k_3\\
&=&S_{1,\sg}(k_1)[S_{1,\sg}(h_1)f(h_2(k_2(x)))h_3]k_3\\
&=&[S_{1,\sg}(k_1)\ced_\sg S_{1,\sg}(h_1)]f((h_2*k_2)(x))(h_3\ced_\sg k_3)\\
&=&[S_{1,\sg}(h_1\ced_{\sg^{-1}}k_1) ]f((h_2*k_2)(x))(h_3\ced_\sg k_3)\\
&=&[f\leftarrow (h k)](x).
\end{array}$$
The third equation holds since $M$ is an $A\#_\sg M$-bimodule. The fourth equation follows from Lemma \ref{anpo}(i). The last equation follows from the fact that both $\bdt_{1,\sg}^\sg$ and $\bdt_{1,1}^\sg$ are algebra homomorphisms. \qed

\begin{rem}\label{rm}
Since $A$ is a graded left $H^\sg$-module algebra,   $A$ is naturally an $A^e\rtimes H^\sg$-module. Hence, $\Hom_{A^e}(A,M)$ is a right $H$-module for any  $(A\#_\sg H)^e$-bimodule $M$. $H$ is just the algebra $H(1,1)$. From the fact that $S_{1,\sg}(h_1)h_2=\vps(h)$ for any $h\in H$, it is easy to check that
$$\Hom_H(\kk,\Hom_{A^e}(A,M))\cong \Hom_{(A\#_\sg H)^e}(A\#_\sg H,M),$$
for any $(A\#_\sg H)^e$-bimodule $M$.
\end{rem}

From Lemma \ref{rh} we see that  $\Hom_{A^e}(P,(A\#_\sg H)^e)$  is a right $H$-module. Moreover, the inner structure of $(A\#_\sg H)^e$ induces a right $(A\#_\sg H)^e$-module structure on $\Hom_{A^e}(P,(A\#_\sg H)^e)$. That is,
$$[f\cdot (a\#h)\ot (b\#g)](x)=f(x)_1(a\#h)\ot(b\#g) f(x)_2$$
for any $f\in \Hom_{A^e}(P,(A\#_\sg H)^e)$ and $a\#h,b\#g\in A\#_\sg H$.

\begin{lem}\label{hh}
Let $P$ be an $A^e\rtimes H^\sg$-module.
\begin{enumerate}\item[(i)]There is a right $H\ot (A\#_\sg H)^e$-module homomorphism
$$\begin{array}{rcl}\Theta:\Hom_{A^e}(P,A^e)_*\ot {}_{*\sg}H\ot {}_\sg H_*&\ra& \Hom_{A^e}(P,(A\#_\sg H)^e)\\
f\ot k\ot l&\mapsto&\Theta(f\ot k\ot l)\end{array}$$
where $\Theta(f\ot k\ot l)(x)=f(x)_1\#k\ot l_1f(x)_2\#l_2$ for any $x\in P$. Here $l_1\ot l_2=\Delta^\sg_{\sg,1}(l)$.
\item[(ii)] If $P$ is finitely generated projective when viewed as an $A^e$-module, then $\Theta$ is an isomorphism.
\end{enumerate}
\end{lem}

In \cite{st}, Stefan showed the relation between the Hochschild cohomologies of A and B, where $B/A$ is a Hopf-Galois extension. When $B=A\#_\sg H$ is a cleft extension, we have the following lemma:
\begin{lem}\cite[Theorem 3.3]{st}\label{spe}
Let $H$ be a Hopf algebra, $\sg$ a 2-cocycle on $H$. Let $A$ be a graded $H^\sg$-module algebra and $N$ an $(A\#_\sg H)^e$-bimodule. Then there is a spectral sequence
$$E_2^{p,q}=\Ext_{H^e}^p(H,\Ext_{A^e}^q(A,N))\Longrightarrow \Ext^{p+q}_{(A\#_\sg H)^e}(A\#_\sg H,N)$$
which is natural in $N$. The right $H$-module $\Ext^q_{A^e}(A,N)$ is viewed as $H^e$-module via the trivial action on the left side.
\end{lem}

\begin{lem}\label{smo}
Let $H$ be a Hopf algebra, $\sg$ a 2-cocycle on $H$ and $A$ a left $H^\sg$-module algebra. If both $A$ and $H$ are homologically smooth, then so is $A\#_\sg H$.
\end{lem}
\proof Let $I$ be an injective $A\#_\sg H$-module. $\Hom_{A^e}(A,I)$  is a right $H$-module by Remark \ref{rm}. From the proof of \cite[Proposition 3.2]{st}, we see that $\Hom_{A^e}(A,I)$ is an injective $H$-module. Moreover, we see in Remark \ref{rm} that $$\Hom_H(\kk,\Hom_{A^e}(A,M))\cong \Hom_{(A\#_\sg H)^e}(A\#_\sg H,M)$$ for any $A\#_\sg H$-bimodule $M$. Now the proof of Proposition 2.11 in \cite{liwz} is valid for the cleft extension $A\#_\sg H$. We obtain that  $A\#_\sg H$ is homologically smooth. \qed

The following lemma is probably well-known, for the convenience of the reader, we provide a proof here.
\begin{lem}\label{kk}
Let $H$ be an augmented algebra such that $H$ is a twisted CY algebra of
dimension $d$ with Nakayama automorphism $\nu$. Then $H$ is of global dimension $d$. Moreover, there is an isomorphism of right $H$-modules
$$\Ext^i_H({}_H\kk,{}_H H)\cong \begin{cases}0,&i\neq d;\\
k_\xi,&i=d,\end{cases}$$
where $\xi:H\ra \kk$ is the homomorphism defined by $\xi(h)=\vps(\nu(h))$ for any $h\in H$.
\end{lem}

\proof If $H$ is an augmented algebra, then $_H\kk$ is a finite
dimensional module. By \cite[Remark 2.8]{bt}, $H$ has global
dimension $d$.

It follows from  \cite[Proposition 2.2]{bt} that $H$ admits a projective
bimodule resolution
$$0\ra P_d\ra \cdots\ra P_1\ra P_0\ra H\ra 0,$$ where each $P_i$ is
finitely generated as an $H$-$H$-bimodule. Tensoring with functor
$\ot_H \kk$, we obtain a projective resolution of
 $_H\kk$:
$$0\ra
P_d\ot_H\kk\ra \cdots\ra P_1\ot_H\kk\ra P_0\ot_H\kk\ra {}_H\kk\ra
0.$$ Since each $P_i$ is finitely generated, the following isomorphisms of right $H$-modules holds:
$$\kk\ot _H\Hom_{H^e}(P_i,H^e)\cong \Hom_{H }(P_i\ot_H\kk,H).$$  Therefore, the complex $\Hom_H(P_\bullet\ot_H\kk,
H)$ is isomorphic to the complex
$\kk\ot_H\Hom_{H^e}(P_\bullet,H^e)$. The algebra $H$ is twisted CY with Nakayama automorphism $\nu$. So the following $H$-$H$-bimodule complex is exact,
$$0\ra \Hom_{H^e}(P_0,H^e)\ra \cdots\ra \Hom_{H^e}(P_{d-1},H^e)\ra  \Hom_{H^e}(P_d,H^e)\ra H^\nu\ra 0.$$ Thus the complex $\kk\ot_H\Hom_{H^e}(P_\bullet,H^e)$ is exact except at $\kk\ot_H\Hom_{H^e}(P_d,H^e)$, whose homology is $\kk\ot_H H^\nu$. It is easy to see that $\kk\ot_H H^\nu\cong \kk_\xi$, where $\xi:H\ra \kk$ is the algebra homomorphism defined by $\xi(h)=\vps(\nu(h))$ for any $h\in H$.
In conclusion, we obtain the following isomorphisms right $H$-modules
$$\Ext_H^i({}_H\kk,{}_HH)\cong\begin{cases}0,&i\neq d;\\\kk_\xi,&i=
d.\end{cases}$$  \qed

\begin{rem}\label{hin} In a similar way, we can also obtain the following isomorphisms of left $H$-modules:
$$\Ext^i_H(\kk_H, H_H)\cong \begin{cases}0,&i\neq d;\\
{}_\eta k,&i=d,\end{cases}$$
where $\eta:H\ra \kk$ is the homomorphism defined by $\eta=\vps\circ\nu^{-1}$.
Therefore, if $H$ is a twisted CY augmented algebra, then $H$ has finite global dimension and satisfy the AS-Gorenstein condition. However, $H$ is not necessarily Noetherian. It is not AS-regular in the sense of Definition \ref{defn as}. We still call $\Ext^i_H({}_H\kk,{}_H H)$ and $\Ext^d_H(\kk_H, H_H)$ left and right homological integral of $H$ and denoted them by $\int^l_H$ and $\int^r_H$ respectively.
\end{rem}

\begin{lem}\label{kk1}
Let $H$ be a twisted CY Hopf algebra with homological integral $\int^l_H=\kk_\xi$, where $\xi:H\ra \kk$ is an algebra homomorphism. Then the Nakayama automorphism $\nu$ of $H$ is given by $\nu(h)=\xi(h_1)S^2(h_2)$ for any $h\in H$. If the right homological integral of $H$ is $\int^l_H={}_\eta \kk$, then $\eta=\xi \circ S$.
\end{lem}
\proof Proposition 4.5(a) in \cite{bz} holds true when the Hopf algebra is not necessarily Noetherian. So we obtain that the Nakayama automorphism $\nu$ satisfies $\nu(h)=\xi(h_1)S^2(h_2)$ for any $h\in H$. From Remark \ref{hin}, we see that $\eta=\vps\circ \nu^{-1}$. Note that for every $h\in H$, $\nu^{-1}(h)=\xi(Sh_1)S^{-2}(h_2)$ and  $\xi\circ S^2(h)=\xi(h)$. Therefore, we obtain that $\eta=\xi \circ S$. \qed

\begin{thm}\label{main}
Let $H$ be a twisted CY Hopf algebra with homological integral $\int^l_H=\kk_\xi$, where $\xi:H\ra \kk$ is an algebra homomorphism and let $\sg$ be a 2-cocycle on $H$. Let $A$ be an $N$-Koszul graded twisted CY algebra with Nakayama automorphism $\mu$ such that $A$ is a left graded $H^\sg$-module algebra. Then $A\#_\sg H$ is a graded twisted CY algebra with Nakayama automorphism  $\rho$ defined by $$\rho(a\#h)=\mu(a)\#\hdet_{H^\sg}(h_1)(S_{\sg,1}^{-1}(S_{1,\sg}^{-1}(h_2)))\xi(h_3)$$ for all $a\#h\in A\#_\sg H$.
\end{thm}
\proof Assume that the CY dimensions of $H$ and $A$ are $d_1$ and $d_2$ respectively.
Take the Koszul complex $\mathbf{K_b(A)}\ra A\ra 0$. In the proof of Proposition \ref{ext}, we see that $\mathbf{K_b(A)}\ra A\ra 0$ is a complex of $A^e\rtimes H^\sg$-modules. It follows from  Lemma \ref{vh1} and Lemma \ref{hh} that the following isomorphisms of $H\ot (A\#H)^e$-module complexes hold:
$$\begin{array}{ccl}\Hom_{A^e}(\mathbf{K_b}(A), (A\#_\sg H)^e)&\cong& \Hom_{A^e}(\mathbf{K_b}(A), (A^e))_*\ot {}_{*\sg }H\ot {}_\sg H_*\\
&\cong &\Hom_{A^e}(\mathbf{K_b}(A), (A^e))\ot {}_\sg H\ot H_*.\end{array}$$
After taking cohomologies, we obtain that
$$\Ext^q_{A^e}(A,(A\#_\sg H)^e)\cong\Ext^q_{A^e}(A,A^e)\ot {}_\sg H\ot H_*$$
as $H\ot (A\#_\sg H)^e$-modules, for any $q\le 0$.

If we view the right $H$-module $\Ext^q_{A^e}(A,(A\#_\sg H)^e)$ as $H^e$-module via the trivial action on the left side, then
$$\begin{array}{ccl}\Ext_{H^e}^p(H,\Ext_{A^e}^q(A,(A\#H)^e))&\cong& \Ext_H^p(\kk,\Ext_{A^e}^q(A,(A\#_\sg H)^e))\\
&\cong& \Ext_H^p(\kk,\Ext^q_{A^e}(A,A^e)\ot {}_\sg H\ot H_*)\\
&\cong& \Ext^q_{A^e}(A,A^e)\ot {}_\sg H\ot\Ext^p_H(\kk_H,H_H).\end{array}$$

By Lemma \ref{spe}, $\Ext^i_{(A\#H)^e}(A\#H,(A\#H)^e))=0$, for $i\neq d_1+d_2$ and
$$
\Ext^{d_1+d_2}_{(A\#H)^e}(A\#H,(A\#H)^e))\cong
\Ext_{A^e}^{d_2}(A,A^e)\ot {}_\sg H \ot\Ext_H^{d_1}(\kk_H, H_H).
 $$
It is an isomorphism of $(A\#_\sg H)^e$-bimodules if the $(A\#_\sg H)^e$-bimodule on $\Ext_{A^e}^{d_2}(A,A^e)\ot {}_\sg H \ot\Ext_H^{d_1}(\kk, H)$ is given by
$$\begin{array}{l}(a\#h)\cdot(x\ot k\ot l)=a((S^{\sg^2} h_1)\rightharpoonup x )\ot S_{1,\sg}(S_{\sg,1}(h_2))\ced_\sg k\ot \xi(Sh_3) l,\\
(x\ot k\ot l)\cdot(b\#g)=x(k_1\cdot b)\ot k_2\ced_\sg g\ot l,\end{array}$$
for any $a\#h,b\#g\in A\#_\sg H$ and $x\ot k\ot l\in \Ext_{A^e}^{d_2}(A,A^e)\ot {}_\sg H \ot\Ext_H^{d_1}(\kk_H, H_H)$.
Note that $\Ext_H^{d_1}(\kk_H, H_H)\cong {}_\eta\kk $, where $\eta=\xi\circ S$ (Lemma \ref{kk1}).

By Proposition \ref{ext}, we obtain the following isomorphism:$$\Ext^{d_1+d_2}_{(A\#H)^e}(A\#H,(A\#H)^e))\cong A_{\mu}\ot A^!_{d_2}\ot H\ot {}_{\xi\circ S}\kk.$$  Since the algebra $A$ is $N$-Koszul graded twisted CY of dimension $d_2$, it is AS-regular of global dimension $d_2$. By \cite[Lemma 5.10]{kkz},  we obtain that $A^!_{d_2}\cong \Ext_A^{d_2}(\kk,\kk)$ is one dimensional. Let $t$ be a nonzero element in $A^!_{d_2}$.  The left $H^\sg$-action on $A^!_{d_2}$ is given by
$$h\cdot t=\hdet(S^{\sg^{-1}}h)t,$$
for any $h\in H$. Therefore, the $(A\#H)^e$-module structure on $A_{\mu}\ot A^!_{d_2}\ot H\ot {}_{\xi\circ S}\kk$ is given by
\begin{equation}
\begin{array}{ll}(a\#h)\cdot(x\ot t\ot k\ot y)&\\
=a(h_1\cdot x)\ot \hdet_{H^\sg}(S^\sg h_2) t\ot (S_{1,\sg} (S_{\sg,1}h_3)) \ced_\sg k\ot \xi(Sh_4)y &\\
(x\ot t\ot k\ot y)\cdot(b\#g)&\\
=x\mu(k_1\cdot b)\ot t\ot k_2\ced_\sg g\ot y, &
\end{array}
\end{equation}
for $(x\ot t\ot k\ot y)\in A_{\mu}\ot A^!_{d_2}\ot H\ot {}_{\xi\circ S}\kk$ and $a\#h,b\#g\in A\#H$.

Now we prove that $A_{\mu}\ot A^!_{d_2}\ot {}_\sg H\ot {}_{\xi\circ S}\kk\cong (A\#_\sg H)^\rho$ as $(A\#_\sg H)^e$-modules for some automorphism $\rho$ of $A\#_\sg H$.

It is straightforward to check that for any $x\in A$, $k\in H$, we have:
$$\begin{array}{ccl}x\ot t\ot k\ot 1 &=& [x\#\hdet_{H^\sg}(k_1)S_{\sg,1}^{-1}(S_{1,\sg}^{-1}(k_2))\xi(k_3)]\cdot(1\ot t\ot 1\ot 1)\\
&=&(1\ot t\ot 1\ot 1)\cdot(\mu^{-1}(x)\#k).
\end{array}$$
This implies that $(1\ot t\ot 1\ot 1)$ is a left and right $A\#_\sg H$-module generator of $A_{\mu}\ot A^!_{d_2}\ot {}_\sg H\ot {}_{\xi\circ S}\kk$. The same formula implies that no nonzero element of $A\#H$ annihilates $(1\ot t\ot 1\ot 1)$. Therefore, $A_{\mu}\ot A^!_{d_2}\ot {}_\sg H\ot {}_{\xi\circ S}\kk$ is a free $A\#_\sg H$-module of rank 1 on each side. So $A_{\mu}\ot A^!_{d_2}\ot {}_\sg H\ot {}_{\xi\circ S}\kk\cong (A\#H)^\rho$ as $(A\#_\sg H)^e$-modules for some automorphism $\rho$ of $A\#_\sg H$. Next we compute $\rho$.  For any $h\in H$,
$$\begin{array}{ccl}(1\ot t\ot 1\ot 1)\cdot(1\#h)&=&1\ot t\ot h\ot 1\\
&=&(1\#\hdet_{H^\sg}(h_1)S_{\sg,1}^{-1}(S_{1,\sg}^{-1}(h_2))\xi(h_3))\cdot(1\ot t\ot 1\ot 1).\end{array}$$
This shows that $\rho(h)=\hdet(h_1)(S_{\sg,1}^{-1}(S_{1,\sg}^{-1}(h_2)))\xi(h_3)$.

On the other hand, for any $a\in A$, we have:
$$\begin{array}{ccl}(1\ot t\ot 1\ot 1)\cdot(a\#1)&=&\mu(a)\ot t\ot 1\ot 1\\
&=&(\mu(a)\#1)\cdot(1\ot t\ot 1\ot 1).\end{array}$$
So $\rho(a)=\mu(a)$. It follows that the automorphism $\rho$ of $A\#_\sg H$ is give by
$$\begin{array}{ccl}\rho(a\#h)&=&\mu(a)\#\hdet_{H^\sg}(h_1)(S_{\sg,1}^{-1}(S_{1,\sg}^{-1}(h_2)))\xi(h_3)
\end{array}$$ for any $a\#h\in A\#H$ and $A_\mu\ot A^!_{d_2}\ot {}_\sg H\ot \kk_\xi\cong { (A\#_\sg H)^\rho}$. To summarize, we obtain the following isomorphisms of $(A\#H)^e$-modules:
$$\Ext^i_{(A\#_\sg H)^e}(A\#_\sg H,(A\#_\sg H)^e)\cong \begin{cases}0, &i\neq d_1+d_2; \\(A\#_\sg H)^\rho,&i=d_1+d_2.\end{cases}$$
By Lemma \ref{smo}, $A\#_\sg H$ is homologically smooth. The proof is completed.
 \qed

Let $H$ be a Hopf algebra.  For an algebra
homomorphism $\xi:H\ra\kk$, We write $[\xi]^l$ for the \textit{left winding homomorphism} of
$\xi$ defined by
$$[\xi]^l(h)=\xi(h_1)h_2,$$ for any $h\in H$. The \textit{right winding automorphism} $[\xi]^r$ of $\xi$ can be defined similarly. It is well-known that both $[\xi]^l$ and $[\xi]^r$ are algebra automorphisms of $H$. In Theorem \ref{main}, if we take the 2-cocycle to be trivial, we obtain the following result about smash products.

\begin{thm}\label{cysma}
Let $H$ be a twisted CY Hopf algebra with homological integral $\int^l_H=\kk_\xi$, where $\xi:H\ra \kk$ is an algebra homomorphism and $A$ an  $N$-Koszul graded twisted CY algebra with Nakayama automorphism $\mu$ such that $A$ is a left graded $H$-module algebra. Then $A\#H$ is a twisted CY algebra with Nakayama automorphism  $\rho=\mu\#(S^{-2}\circ[\hdet_H]^l\circ [\xi]^r)$.
\end{thm}
\proof From Theorem \ref{main}, we see that  $A\#H$ is a graded twisted CY algebra with Nakayama automorphism  $\rho$ defined by $$\rho(a\#h)=\mu(a)\#\hdet_H(h_1)(S^{-2}(h_2))\xi(h_3)$$ for all $a\#h\in A\#_\sg H$. That is, $\rho=\mu\#(S^{-2}\circ[\hdet_H]^l\circ [\xi]^r)$. \qed

\begin{cor}\label{cor}
With the same assumption as in Theorem \ref{cysma}, the algebra $A\#H$ is a CY algebra if and only if $\hdet_H=\xi\circ S$ and $\mu\#S^{-2}$ is an inner automorphism of $A\#H$.
\end{cor}
\proof Since $\mu\#(S^{-2}\circ [\hdet_H]^l\circ [\xi]^r)=(\mu\#S^{-2})\circ (\id\#([\hdet_H]^l\circ [\xi]^r))$, the sufficiency part is clear.

In the proof of Theorem \ref{main}, if we let the cocycle $\sg$  be trivial, then the proof is just a modification of the proof of the sufficiency part of \cite[Theorem 2.12]{liwz}. If we modify the proof of the necessary part,   we obtain that  $\xi\star \hdet_H=\vps$, where $\star$ stands for the convolution product.  It is easy to see that $\xi\circ S$ and $\xi$ are inverse to each other with respect to the convolution product. Therefore, we obtain that $\hdet_H =\xi \circ S$. Now $\mu\#(S^{-2}\circ [\hdet_H]^l\circ [\xi]^r)=\mu\#S^{-2}$. It follows from  Theorem \ref{cysma} that $\mu\#S^{-2}$ is an inner automorphism.  \qed

In case $A$ is an $N$-Koszul graded CY algebra and $H$ is a CY Hopf algebra,  we have the following consequence.

\begin{cor}
Let $H$ be a CY Hopf algebra, and let $A$ be an $N$-Koszul graded CY algebra and a left graded $H$-module algebra. Then $A\#H$ is a graded CY algebra if and only if the homological determinant of the $H$-action on $A$ is trivial and $\id\#S^2$ is an inner automorphism of $A\#H$.
\end{cor}
\proof Since $H$ is a CY Hopf algebra, by Lemma \ref{cyas} (ii), the algebra $H$ satisfies $\int_H^l=\kk$. Now the corollary  follows immediately from Corollary \ref{cor}.  \qed

\begin{rem}
From Lemma \ref{kk} and Lemma \ref{kk1}, it is not hard to see that if $H$ is CY Hopf algebra, then $S^2$ is an inner automorphism of $H$. However, $\id\#S^2$ is not necessarily an inner automorphism of $A\#H$ even if $A\#H$ is CY. Example \ref{eg1} in Section \ref{s4} is  a counterexample. It also shows that the smash product $A\#H$ could be a CY algebra when $A$ itself is not.
\end{rem}

In Theorem \ref{main}, if we let the algebra $A$ be $\kk$, then we obtain the following result about the twisted CY property of cleft objects.
\begin{thm}\label{cleft}
Let $H$ be a twisted CY Hopf algebra with $\int^l_H={}_\xi\kk$.  Suppose ${}_\sg H$ is a right cleft object of $H$. Then ${}_\sg H$ is a twisted CY algebra with Nakayama automorphism $\mu$ defined by
$$\mu(x)=S^{-1}_{\sg,1}(S^{-1}_{1,\sg}(x_1))\xi S(x_2)$$
for any $x\in {}_\sg H$.
\end{thm}

\section{Cleft objects of $U(\mc{D},\lmd)$}\label{s3}
The pointed Hopf algebras $U(\mc{D},\lmd)$ introduced in \cite{as3} are generalizations of the quantized enveloping algebras $U_q(\gg)$, where $\gg$ is a finite dimensional semisimple Lie algebra. Chelma showed that the algebras  $U_q(\gg)$ are CY algebras \cite[Theorem 3.3.2]{c}.
The CY property of the algebras $U(\mc{D},\lmd)$ were discussed in \cite{yz}. In this section we will show that the  cleft objects of the algebras $U(\mc{D},\lmd)$ are all twisted CY algebras.

\subsection{The Hopf algebra $U(\mc{D},\lmd)$} We refer to \cite{as} for a detailed discussion about braided Hopf algebras and Yetter-Drinfeld modules. For a group $\bgm$, we denote by  ${}_\bgm^\bgm \mathcal {YD}$
the category of Yetter-Drinfeld modules over the group algebra $\kk\bgm$. If $\bgm$ is an abelian group, then it is well-known that a Yetter-Drinfeld module over the algebra $\kk\bgm$ is just a $\bgm$-graded $\bgm$-module.

%For $M\in {}_\bgm^\bgm \mathcal {YD}$, $g\in \bgm$ and $\chi\in \hat{\bgm}$

 We fix the following terminology.
\begin{itemize}
\item  a free abelian group $\bgm$ of finite rank $s$;
\item  a Cartan matrix $\mathbb{A}=(a_{ij})\in \mathbb{Z}^{\tt\times \tt}$ of finite type, where $\tt\in\NN$. Let  $(d_1,\cdots, d_\tt)$ be a diagonal matrix
of positive integers such that $d_ia_{ij} = d_ja_{ji}$, which is
minimal with this property;
\item a set $\mathcal {X}$ of connected components of the Dynkin diagram corresponding
to the Cartan matrix $\mathbb{A}$. If $1\se i, j\se \tt$, then $i\sim
j$ means that they belong to the same connected component;
\item a family $(q_{_I})_{I\in \mc{X}}$ of elements in $\kk$ which are \textit{not} roots of unity;
\item elements $g_1,\cdots , g_\tt\in \bgm$ and characters $\chi_1,\cdots, \chi_\tt\in \hat{\bgm}$ such that
\begin{equation}\label{q}\chi_j(g_i)\chi_i(g_j)=q_{I}^{d_ia_{ij}}, \t{  } \chi_i(g_i)=q_{I}^{d_i}, \t{   for all $1\se i,j\se \tt$, $I\in
\mc{X}$}.\end{equation}
\end{itemize}

For simplicity, we write $q_{ji}=\chi_i(g_j)$. Then Equation (\ref{q}) reads as follows:
\begin{equation}\label{q1}q_{ii}=q_I^{d_i}\t{ and } q_{ij}q_{ji}=q_{I}^{d_ia_{ij}}\t{ for all }
1\se i,j\se \tt, I\in \mc{X}.
\end{equation}

Let $\mc{D}$ be the collection $\mc{D}(\bgm,(a_{ij})_{1\se i,j\se
\tt}, (q_{_{I}})_{I\in \mathcal {X}}, (g_i)_{1\se i\se
\tt},(\chi_i)_{1\se i\se \tt} )$. A \textit{linking datum}
$\lmd=(\lmd_{ij})$ for $\mc{D}$ is a collection of elements
$(\lmd_{ij})_{1\se i<j\se \tt,i\nsim j}\in\kk$ such that
$\lmd_{ij}=0$ if $g_ig_j=1$ or $\chi_i\chi_j\neq\varepsilon$. We
write the datum $\lmd=0$, if $\lmd_{ij}=0$ for all $1\se i<j\se
\tt$. The datum $(\mc{D},\lmd)=(\bgm,(a_{ij}) , q_{_I},
(g_i),(\chi_i), (\lmd_{ij}) )$ is called a \textit{ generic datum of
finite Cartan type} for group $\bgm$.

A generic datum of finite Cartan type for a group $\bgm$ defines a Yetter-Drinfeld module over the group algebra $\kk\bgm$. Let $V$ be a vector space with basis $\{x_1,x_2,\cdots,x_\tt\}$. We set
$$|x_i|=g_i,\;\;g(x_i)=\chi_i(g)x_i, \;\;1\se i\se \tt, g\in \bgm,$$
where $|x_i|$ denote the degree of $x_i$.
This makes $V$ a Yetter-Drinfeld module over the group algebra $\kk\bgm$. We write $V=\{x_i,g_i,\chi_i\}_{1\se i\se \tt}\in {}^\bgm_\bgm\mc{YD}$.
The braiding is given
by
$$c(x_i\ot x_j)=q_{ij}x_j\ot x_i, \;\;1\se i,j\se \tt.$$

The tensor algebra $T(V)$ on $V$ is a natural graded braided Hopf algebra in ${}_\bgm^\bgm \mathcal {YD}$. The smash product $T(V)\#\kk\bgm$ is a usual Hopf algebra. It is also called a bosonization of $T(V)$ by $\kk\bgm$.

\begin{defn} Given a generic datum of finite Cartan type $(\mc{D},\lmd)$ for a group $\bgm$. Define $U(\mc{D},\lmd)$ as the quotient Hopf algebra of the smash product $T(V)\#\kk\bgm$ modulo the ideal generated by
$$(\t{ad}_cx_i)^{1-a_{ij}}(x_j)=0,\;\;1\se i\neq j\se \tt,\;\;i\sim j,$$
$$x_ix_j-\chi_j(g_i)x_jx_i=\lmd_{ij}(g_ig_j-1),\;\;1\se
i<j\se \tt,\;\;i\nsim j,$$
where $\t{ad}_c$ is the braided adjoint representation defined  in
\cite[Sec. 1]{as3}.
\end{defn}

The algebra $U(\mc{D},\lmd)$ is a pointed Hopf algebra with
$$\bdt(g)=g\ot g,\; \bdt(x_i)=x_i\ot 1+g_i\ot x_i,\;\; g\in \bgm, 1\se i\se \tt.$$

To present the CY property of the algebras $U(\mc{D},\lmd)$, we recall the concept of root vectors.
Let $\Phi$ be the root system corresponding to the Cartan matrix
$\mathbb{A}$ with $\{\al_1,\cdots, \al_\tt\}$ a set of fix simple
roots, and  $\mathcal {W}$ the Weyl group. We fix a reduced
decomposition of the longest element $w_0=s_{i_1}\cdots s_{i_p}$ of
$\mathcal {W}$ in terms of the simple reflections.  Then the
positive roots are precisely the followings,
$$\bt_1=\al_{i_1}, \;\;\bt_2=s_{i_1}(\al_{i_2}),\cdots, \bt_p=s_{i_1}\cdots s_{i_{p-1}}(\al_{i_p}).$$ For
$\bt_i=\sum_{i=1}^{\tt} m_i\al_i$, we write
$$g_{\bt_i}=g_1^{m_1}\cdots g_\tt^{m_\tt} \t{ and }\chi_{\bt_i}={\chi}_1^{m_1}\cdots {\chi}_\tt^{m_\tt}.$$

Lusztig defined the root vectors for a quantum group $U_q(\mathfrak{g})$ in
 \cite{l}. Up to a non-zero scalar, each root vector can be
expressed as an iterated braided commutator. In \cite[Sec.
4.1]{as4}, the root vectors were generalized on a pointed Hopf
algebras $U(\mc{D},\lmd)$. For each positive root $\bt_i$, $1\se i
\se p$, the root vector $x_{\bt_i}$ is defined by the same iterated
braided commutator of the elements $x_1, \cdots , x_\tt$, but with
respect to the general braiding.

\begin{rem}\label{root} If $\bt_j=\al_l$, then we have $x_{\bt_j}=x_l$. That is, $x_1,\cdots,x_\tt$ are the simple root vectors.
\end{rem}

\begin{lem}\label{cy}
Let $(\mc{D},\lmd)$ be a generic datum of finite Cartan type for a group $\bgm$, and $H$
the Hopf algebra $U(\mc{D},\lmd)$. Let $s$ be the rank of $\bgm$ and $p$
the number of the positive roots of the Cartan matrix.
\begin{enumerate}
\item[(i)] The algebra $H$ is Noetherian AS-regular
of global dimension $p+s$. The left
homological integral module $\int^l_H$ of $H$ is isomorphic to
$\kk_\zeta$, where $\zeta:H\ra \kk$ is an algebra homomorphism defined
by $\zeta(g)=(\prod_{i=1}^p\chi_{_{\bt_i}})(g)$ for all $g\in \bgm$
and $\zeta(x_k)=0$ for all $1\se k\se \tt$.
\item[(ii)] The algebra $H$ is twisted CY with Nakayama automorphism $\mu$ defined by  $
\mu(x_k)=q_{kk}x_k$, for all
$1\se k\se \tt$, and $\mu(g)=(\prod_{i=1}^p\chi_{_{\bt_i}})(g)$ for
all $g\in \bgm$.
\item[(iii)]  The algebra $H$ is CY if and only if
$\prod_{i=1}^p\chi_{_{\bt_i}}=\varepsilon$ and $S^2$ is an
inner automorphism.
\end{enumerate}
\end{lem}
\proof (i) This is Theorem 2.2 in \cite{yz}.

(ii) By Lemma \ref{cyas}(i),  we conclude that the algebra $H$ is twisted CY with Nakayama automorphism $\mu$ defined
by $\mu(x_k)=S^{-2}(x_k)=q_{kk}x_k$ for $1\se k\se \tt$ and  $\mu(g)=\xi(g)g=(\prod_{i=1}^p\chi_{_{\bt_i}})(g)g$ for $g\in \bgm$.

(iii) This follows directly from (i) and Lemma \ref{cyas} (ii). \qed

\begin{rem}
Theorem 2.3 in \cite{yz} showed that the Nakayama automorphism of the algebra $U(\mc{D},\lmd)$ is the algebra automorphism  $\nu$ defined by $\nu(x_k)=\prod_{i=1,i\neq j_k}^{p}\chi_{_{\bt_i}}(g_k)x_k$, for all
$1\se k\se \tt$, and $\nu(g)=(\prod_{i=1}^p\chi_{_{\bt_i}})(g)$ for
all $g\in \bgm$, where each $ j_k $, $1\se k\se \tt$,  is the integer
such that $\bt_{j_k}=\al_k$. Now we show that the algebra automorphisms $\mu$ and $\nu$ only differ by an inner automorphism.

By a similar discussion to the one in the proof of Lemma 4.1 in \cite{yz}, we see that
$$\prod_{i=1,i\neq j_k}^{p}\chi_{_{\bt_i}}(g_k)=(\prod_{i=1}^{j_k-1}\chi^{-1}_k(g_{_{\bt_i}}))( \prod_{i=j_k+1}^{p}\chi_{_{\bt_i}}(g_k))=\prod_{i=1,i\neq j_k}^{p}\chi^{-1}_k(g_{\bt_i})$$
for each $1\se k\se \tt$.
Therefore, $$\begin{array}{ccl}[\prod_{i=1}^{p}g_{\bt_i}]^{-1}(\mu(x_k))[\prod_{i=1}^{p}g_{\bt_i}]&=&\prod_{i=1}^{p}\chi_k^{-1}(g_{\bt_i})q_{kk}x_k\\
&=&\prod_{i=1,i\neq j_k}^{p}\chi^{-1}_k(g_{\bt_i})x_k\\&=&\prod_{i=1,i\neq j_k}^{p}\chi_{_{\bt_i}}(g_k)\\&=&\nu(x_k)\end{array}$$ for $1\se k\se \tt$. Moreover, $\bgm$ is abelian, so $[\prod_{i=1}^{p}g_{\bt_i}]^{-1}(\mu(g))[\prod_{i=1}^{p}g_{\bt_i}]=\mu(g)=\nu(g)$ for all $g\in \bgm$. This shows that  $\mu$ and $\nu$ indeed differ by an inner automorphism.
\end{rem}

In \cite{ma}, the author classified the cleft objects of a class of pointed Hopf algebras. This class of algebras contains the algebras $U(\mc{D}, \lmd)$.

Now we fix  a generic  datum of finite Cartan
type
$$(\mc{D},\lmd)=(\bgm,(a_{ij})_{1\se i,j\se \tt},(q_{_I})_{I\in \mathcal {X}},
(g_i)_{1\se i\se \tt},(\chi_i)_{1\se i\se \tt},(\lmd_{ij})_{1\se
i<j\se \tt,i\nsim j} ),$$  where $\bgm$  is a free abelian group of
rank $s$.

Let $\sg\in Z^2(\kk\bgm)$ be a 2-cocycle for the group algebra $\kk\bgm$. Define $\chi_i^\sg(g)=\frac{\sg(g,g_i)}{\sg(g_i,g)}\chi_i(g)$. From \cite[Proposition 1.11]{ma}, we obtain that $${}_\sg V=\{x_i,g_i,\chi_i^\sg\}_{1\se i\se\tt}\in {}^\bgm_\bgm \mc{YD}.$$
The associated braiding  is given by
$$c^\sg(x_i\ot x_j)=q^\sg_{ij}x_j\ot x_i,$$
where $q^\sg_{ij}=\frac{\sg(g_i,g_j)}{\sg(g_j,g_i)}q_{ij}$.

Define $$\Xi(\sg)=\{(i,j)\;|\; i<j,i\nsim j,\chi^\sg_i\chi^\sg_j=1\}.$$

Given the braided vector space ${}_\sg V$, we have the tensor algebra $T({}_\sg V)$ and the smash product $T({}_\sg V)\#\kk\bgm$. The  2-cocycle $\sg$ for the group algebra $\kk\bgm$ can be regarded as a 2-cocycle for $T({}_\sg V)\#\kk\bgm$ through the projection $T({}_\sg V)\#\kk\bgm\ra \kk\bgm$. Then we have the crossed product $T({}_\sg V)\#_\sg\kk\bgm$. The difference between the crossed product and the smash product $T({}_\sg V)\#\kk\bgm$ is given by
$$\ol{g}\ol{g'}=\sg(g,g')\ol{gg'}, \;\;g,g'\in\bgm, \forall g\in G. $$
Here $g\in T({}_\sg V)\#\kk\bgm$ is denoted by $\ol{g}\in T({}_\sg V)\#_\sg\kk\bgm$ to avoid confusion.

\begin{defn}
Given  $\pi=(\pi_{ij})\in \kk^{\Xi(\sg)}$. Define $B^\lmd(\sg,\pi)$ to be the quotient algebra of $T({}_\sg V)\#_\sg\kk\bgm$ modulo the ideal generated by
$$(\ad_{c^\sg} x_i)^{1-a_{ij}}(x_j)=0,\;\; 1\se i\neq j\se \tt,\;\;i\sim j,$$
$$(\ad_{c^\sg} x_i)(x_j)-\lmd_{ij}\ol{g}_i\ol{g}_j+\pi_{ij}=0, \;\;1\se i<j\se \tt, i\nsim j,$$
where we set $\pi_{ij}=0$ if $(i,j)\notin\Xi(\sg) $.
\end{defn}

Let $\mc{Z}=\mc{Z}(\bgm,\Xi,\kk)$ denote the set of all pairs $(\sg,\pi)$, where $\sg\in Z^2(\kk\bgm)$ and $\pi=(\pi_{ij})\in \kk^{\Xi(\sg)}$.
For two pairs $(\sg, \pi)$ and $(\sg',\pi')$, define $(\sg, \pi)\thicksim(\sg',\pi')$, if there is an invertible map $f:\kk\bgm\ra \kk$ such that
$$\sg'(g,h)=f^{-1}(g)f^{-1}(h)\sg(g,h)f(gh), \;\; g,h\in \bgm;$$
$$\pi'_{ij}=f^{-1}(g_i)f^{-1}(g_j)\pi_{ij},\;\;(i,j)\in \Xi(\sg).$$ This defines an equivalence relation on $\mc{Z}$. We write
$\mc{H}(\bgm,\Xi,\kk)=\mc{Z}/\thicksim$.

The following Lemma is the right version of Theorem 6.3 in \cite{ma}. It describes the isomorphism classes of right cleft objects of the algebras $U(\mc{D},\lmd)$.
\begin{lem}\label{cle}
The map defined by
$$\begin{array}{rcl}
\mc{H}(\bgm,\Xi,\kk)&\lra &\Cleft(U(\mc{D},\lmd))\\
(\sg,\pi)&\longmapsto&B^\lmd(\sg,\pi)
\end{array}$$
is a bijection, where $\Cleft(U(\mc{D},\lmd))$ denotes the set of the isomorphism classes the right cleft objects of $U(\mc{D},\lmd)$.
\end{lem}

\begin{prop}\label{cycle}
Given a pair $(\sg, \pi)\in \mc {Z}(\bgm, \Xi,\kk)$. The algebra $B^{\lmd}(\sg,\pi)$ is twisted CY with Nakayama automorphism defined by
$\mu(x_k)=q_{kk}x_k$ for all $1\se k\se \tt$ and $\mu(g)=(\prod_{i=1}^p\chi_{_{\bt_i}})(g)$ for
all $g\in \bgm$.

In particular, the algebra $B^{\lmd}(\sg,\pi)$ is CY if and only if there is en element $h\in\kk\bgm$ such that
$\frac{\sg(h,g)}{\sg(g,h)}=(\prod_{i=1}^p\chi_{_{\bt_i}})(g)$, for all $g\in \bgm$ and $(\prod_{i=1,i\neq j_k}^p\chi_{_{\bt_i}})(g)\chi_k(h)=1$
for each $1\se k\se \tt$, where each $ j_k $, $1\se k\se \tt$,  is the integer
such that $\bt_{j_k}=\al_k$.
\end{prop}
\proof Let $H=U(\mc{D},\lmd)$. Without loss of generality, we may assume that $\sg$ satisfies that
$$\sg(g,g^{-1})=\sg(g^{-1},g)=1$$ for all $g\in \bgm$.
This follows from Lemma \ref{cle} and the fact that for
 each pair $(\sg,\pi)$, there is a pair $(\sg',\pi')$ such that $(\sg,\pi)\sim (\sg',\pi')$ and $\sg'$ satisfies
$\sg'(g,g^{-1})=\sg'(g^{-1},g)=1$ for all $g\in \bgm$.  The algebra $B^\lmd_q(\sg,\pi)$ is a cleft object of $H$. Then $B^\lmd_q(\sg,\pi)\cong {}_\tau H$, for some 2-cocycle $\tau$. The 2-cocycle $\tau$ can be calculated using Lemma \ref{mm}. We conclude that  $\tau$ satisfies the following:
$$\begin{array}{rcl}\tau(g,g')&=&\sg(g,g'),\\\tau(g,x_i)&=&\tau(x_i,g)=0,\;\;1\se i\se \tt, g,g'\in \bgm.\\\tau(x_i,x_j)&=&\begin{cases}\lmd_{ij}\sg(g_i,g_j)-\pi_{ij},&i<j,i\nsim j\\0,&otherwise.\end{cases}\end{array}$$

Lemma \ref{cy} shows that the algebra $H=U(\mc{D},\lmd)$ is Noetherian AS-regular. The left
homological integral module $\int^l_H$ of $H$ is isomorphic to
$\kk_\zeta$, where $\zeta:H\ra \kk$ is an algebra homomorphism defined
by $\zeta(g)=(\prod_{i=1}^p\chi_{_{\bt_i}})(g)$ for all $g\in \bgm$
and $\zeta(x_k)=0$ for all $1\se k\se \tt$.

Since $H$ is AS-regular, by Theorem \ref{cleft}, $B_q(\sg,\pi)\cong {}_\tau H$ is a twisted CY algebra. Its Nakayama automorphism can be calculated as follows.
For $g\in \bgm$, $$\begin{array}{ccl}\mu(g)&=&S^{-1}_{\tau,1}(S^{-1}_{1,\tau}(g))\zeta(g)=S^{-1}_{\tau,1}(g^{-1}\sg(g^{-1},g))\zeta(g)\\
&=&S^{-1}_{\tau,1}(g^{-1})\zeta(g)=\sg(g,g^{-1})g\zeta(g)=\zeta(g)g\\
&=&(\prod_{i=1}^p\chi_{_{\bt_i}})(g)g.\end{array}$$
For each $1\se k\se \tt$,
$$\begin{array}{ccl}\mu({x}_k)&=&S^{-1}_{\tau,1}(S^{-1}_{1,\tau}({x}_k))=S^{-1}_{\tau,1}(-g^{-1}_k{x}_k\sg(g_k^{-1},g_k))\\
&=&S^{-1}_{\tau,1}(-g_k^{-1}{x}_k)=\sg(g_k^{-1},g_k)q_{kk}x_k\\&=&q_{kk}x_k.\end{array}$$

The algebra $B^{\lmd}(\sg,\pi)$ is CY if and only if the algebra automorphism $\mu$ is inner. Since   the algebra $U(\mc{D},\lmd)$ is a domain \cite[Theorem 4.3]{as3}, the invertible elements of $B^{\lmd}(\sg,\pi)$ fall in $\kk\bgm$. In $B^{\lmd}(\sg,\pi)$, for $l,g\in \bgm$ and $1\se k\se \tt$, we have
$$\ol{l}\ol{g}=\frac{\sg(l,g)}{\sg(g,l)}\ol{g}\ol{l},\;\;\ol{l}x_k=\chi^{\sg}_k(l)x_k \ol{l}=\frac{\sg(l,g_k)}{\sg(g_k,l)}\chi_k(l)x_k \ol{l}.$$
With these facts, we see that the automorphism $\mu$ is an inner automorphism if and only if there exists an element $h\in \kk\bgm$ such that
\begin{equation}\label{2}\frac{\sg(h,g)}{\sg(g,h)}=(\prod_{i=1}^p\chi_{_{\bt_i}})(g),\;\;\frac{\sg(h,g_k)}{\sg(g_k,h)}\chi_k(h)=q_{kk},\end{equation}
for all $g\in \bgm$ and $1\se k\se \tt$. Note that if $\frac{\sg(h,g)}{\sg(g,h)}=(\prod_{i=1}^p\chi_{_{\bt_i}})(g)$ holds for any $g\in \bgm$, then $\frac{\sg(h,g_k)}{\sg(g_k,h)}=(\prod_{i=1}^p\chi_{_{\bt_i}})(g)$. So the condition (\ref{2}) is equivalent to
$$\frac{\sg(h,g)}{\sg(g,h)}=(\prod_{i=1}^p\chi_{_{\bt_i}})(g),(\prod_{i=1,i\neq j_k}^p\chi_{_{\bt_i}})(g)\chi_k(h)=1,$$ for all $g\in \bgm$ and
 $1\se k\se \tt$, where each $ j_k $, $1\se k\se \tt$,  is the integer
such that $\bt_{j_k}=\al_k$.
\qed

We end this section by giving some examples. We first need the following lemma.

\begin{lem}\label{u}
Let $\bgm$ be an abelian group, $\sg$ a 2-cocycle for the group algebra $\kk\bgm$. For any $g,k,h\in \bgm$, we have
$$\frac{\sg(gk,h)}{\sg(h,gk)}=\frac{\sg(g,h)}{\sg(h,g)}\frac{\sg(k,h)}{\sg(h,k)}.$$
\end{lem}
\proof Since $\sg$ is a 2-cocycle, the following equations hold for any $g,h,k\in \bgm$.
\begin{equation}\label{3} \sg(g,k)\sg(gk,h)=\sg(k,h)\sg(g,kh)\end{equation}
\begin{equation}\label{4} \sg(g,k)\sg(h,gk)=\sg(h,g)\sg(hg,k)\end{equation}
\begin{equation}\label{5} \sg(g,h)\sg(gh,k)=\sg(h,k)\sg(g,hk)\end{equation}
\begin{equation}\label{6} \sg(h,k)\sg(g,hk)=\sg(g,h)\sg(gh,k)\end{equation}
By (\ref{3}) and (\ref{4}), we obtain
$$\begin{array}{ccl}\frac{\sg(gk,h)}{\sg(h,gk)}&=&\frac{\sg(k,h)\sg(g,kh)}{\sg(h,g)\sg(hg,k)}\\
&\overset{(\ref{5},\ref{6})}=&\frac{\sg(k,h)\frac{\sg(g,h)\sg(gh,k)}{\sg(h,k)}}{\sg(h,g)\frac{\sg(h,k)\sg(g,hk)}{\sg(g,h)}}\\
&=&\frac{\sg(g,h)}{\sg(h,g)}\frac{\sg(k,h)}{\sg(h,k)}\frac{\sg(g,h)\sg(gh,k)}{\sg(h,k)\sg(g,hk)}\\
&\overset{(\ref{5})}=&\frac{\sg(g,h)}{\sg(h,g)}\frac{\sg(k,h)}{\sg(h,k)}.
\end{array}$$\qed

Now we give an example in which the algebra $U(\mc{D},\lmd)$ is CY, but the algebra $B^\lmd(\sg,\pi)$ is not necessarily CY.

\begin{eg}
 Let  $(\mc{D},\lmd)$  be the datum  given
by \begin{itemize}
\item $\bgm=\lan y_1,y_2\ran$, a free abelian group of rank 2;
\item The Cartan matrix is of type $A_2\times A_2$;
\item $g_1=g_3=y_1,g_2=g_4=y_2$;
\item $\chi_1(y_1)=q^2,\chi_1(y_2)=q^{-1},\chi_2(y_1)=q^{-1},\chi_2(y_2)=q^{-2}$, and $\chi_3=\chi_1^{-1}, \chi_4=\chi_2^{-1}$, where $q$ is not a root of unity;
\item $\lmd=(\lmd_{13},\lmd_{14},\lmd_{23},\lmd_{24})=(0,1,1,0)$.\end{itemize}

Then the algebra $U(\mc{D},\lmd)$ is just the quantized enveloping algebra $U_q(\gg)$, where $\gg$ is the simple Lie algebra corresponding to the Cartan matrix of type $A_2$. Therefore, $U(\mc{D},\lmd)$ is CY (\cite[Theorem 3.3.2]{c}). In fact, we have that
 $$\bt_1=\al_1,\;\;\bt_2=\al_1+\al_2,\;\;\bt_3=\al_2, \;\;\bt_4=\al_3,\;\;\bt_5=\al_3+\al_4,\;\;\bt_6=\al_4$$
are the positive roots, where $\al_i$ $(1\se i\se 4)$ are the simple roots. Hence $\prod_{i=1}^6\chi_{_{\bt_i}}=\chi^2_1\chi^2_2\chi^2_3\chi^2_4=\vps$.
Moreover, $(y^{-2}_1y_2^{-2})x_i(y^2_1y_2^2)=q_{ii}^{-1}x_i=S^2(x_i)$ for $1\se i\se 4$.

Let $\sg$ be a 2-cocycle such that $u_{12}=\frac{\sg(y_2,y_1)}{\sg(y_1,y_2)}$ is not a root of unity. Let $u_{21}=u_{12}^{-1}$. We claim that the algebra $B^\lmd(\sg, \pi)$ can not be a CY algebra. Otherwise, by Proposition \ref{cycle}, there is an element $y_1^iy_2^j\in \bgm$ such that for any $y_1^ky_2^l\in \bgm$, $\frac{\sg(y_1^iy_2^j,y_1^ky_2^l)}{\sg(y_1^ky_2^l,y_1^iy_2^j)}=u_{21}^{il}u_{12}^{jk}=1$, where the first equation follows from Lemma \ref{u} and the second equation holds because $\prod_{i=1}^6\chi_{_{\bt_i}}=\vps$. Now let $k=l=1$. We obtain that $u_{21}^iu_{12}^j=u_{12}^{i-j}=1$, Since $u_{12}$ is not a root of unity, we have that $i=j$. Then $u_{21}^{il}u_{12}^{jk}=u_{12}^{k-l}$ can not equal to 1 when $k\neq l$. This is a contradiction.
\end{eg}

The next example shows  that the algebra $U(\mc{D},\lmd)$ is not CY, but some cleft objects are CY.

\begin{eg}
 Let  $(\mc{D},\lmd)$  be the datum  given
by
\begin{itemize}
\item $\bgm=\lan y_1,y_2\ran$, a free abelian group of rank 2;
\item The Cartan matrix $\mathbb{A}$ is of type $A_1\times A_1$;
\item $g_1=y_1,g_2=y_2$;
\item $\chi_1(g_1)=q^2,\chi_1(g_2)=q^{-4},\chi_2(g_1)=q^4,\chi_2(g_2)=q^{-2}$, where $q$ is not a root of unity;
\item $\lmd=0$.\end{itemize}

The positive roots of $\mathbb{A}$ are just the simple roots. Since $\chi_1\chi_2\neq \vps$, the algebra $H=U(\mc{D},\lmd)$ is not CY (Lemma \ref{cy} (c)).

Let $B^0(\sg,\pi)$ be a cleft object of $H$ such that the 2-cocycle $\sg$ satisfies $u_{12}=\frac{\sg(g_2,g_1)}{\sg(g_1,g_2)}=q^3$. We also put $u_{21}=u^{-1}_{12}$. Choose an element $h=g_1^2g_2^2\in \bgm$ . Then $$\frac{\sg(h,g_1)}{\sg(g_1,h)}=\frac{\sg(g_1^2g_2^2,g_1)}{\sg(g_1,g_1^2g_2^2)}=u_{12}^2=q^6=\chi_1\chi_2(g_1),$$
where the second equation also follows from Lemma \ref{u}. Similarly,
$$\frac{\sg(h,g_2)}{\sg(g_2,h)}=\frac{\sg(g_1^2g_2^2,g_2)}{\sg(g_2,g_1^2g_2^2)}=u_{21}^2=q^{-6}=\chi_1\chi_2(g_2).$$
Moreover, $$\chi_2(g_1)\chi_1(h)=\chi_2(g_1)\chi_1(g_1^2g_2^2)=1,$$
$$\chi_1(g_2)\chi_2(h)=\chi_1(g_2)\chi_2(g_1^2g_2^2)=1.$$
By Proposition \ref{cycle}, the algebra $B^0(\sg,\pi)$ is a CY algebra.
\end{eg}

\section{More Examples}\label{s4}
In this section, we give some examples of Theorem \ref{main}.

The following example shows that it is possible that the crossed product of  CY algebras might be a CY algebra, while their smash product is not CY.

\begin{eg}
Let $A=\kk\lan x_1,x_2\ran/(x_1x_2-x_2x_1)$ be the polynomial algebra with two variables. Then $A$ is a CY algebra. Let $\bgm$ be the free abelian group of rank 2 with generators $g_1$ and $g_2$. There is a $\bgm$-action on $A$ as follows:
$$\begin{array}{cccc}
g_1\cdot x_1=q x_1,&g_2\cdot x_1=q^{-1} x_1,\\
g_1\cdot x_2=qx_2,&g_2\cdot x_2=q^{-1}x_2,
\end{array}$$
where $q$ is not a  root of unity. The homological determinant of this $\bgm$-action is not trivial, namely, $\hdet(g_1)=q^2$,  $\hdet(g_2)=q^{-2}$. The algebra $A\# \kk \bgm$ is not a CY algebra by Theorem 2.12 in \cite{liwz}.

Let $\sg$ be a 2-cocycle on $\bgm$ such that $\frac{\sg(g_2,g_1)}{\sg(g_1,g_2)}=q$. Without loss of generality, we may assume that $\sg(g,g^{-1})=\sg(g^{-1},g)=1$ for $g\in \bgm$. Then the algebra $A\#_\sg \kk\bgm$ is a twisted CY algebra with Nakayama automorphism $\rho$ defined by $\rho(a\#g)=\hdet(h)a\#g$ for any $a\#g\in A\#_\sg \kk\bgm $. Choose an element $h=g_1^2g_2^2\in \bgm$. By Lemma \ref{u}, $$\frac{\sg(h,g_1)}{\sg(g_1,h)}=\frac{\sg(g_1^2g_2^2,g_1)}{\sg(g_1,g_1^2g_2^2)}=(\frac{\sg(g_2,g_1)}{\sg(g_1,g_2)})^2=q^2=\hdet(g_1),$$
$$\frac{\sg(h,g_2)}{\sg(g_2,h)}=\frac{\sg(g_1^2g_2^2,g_2)}{\sg(g_2,g_1^2g_2^2)}=(\frac{\sg(g_1,g_2)}{\sg(g_2,g_1)})^2=q^{-2}=\hdet(g_2).$$
Moreover, $h\cdot x_i=x_i$, $1\se i\se 2$. Therefore, $\rho(a\#g)=h(a\#g)h^{-1}$, for any $a\#g\in A\#_\sg \kk\bgm $. The automorphism $\rho$ is an inner automorphism. So the algebra $A\#_\sg \kk\bgm $ is a CY algebra.
\end{eg}

In the followings, we provide some examples involving the algebras $U(\mc{D},\lmd)$. The definitions of algebras $U(\mc{D},\lmd)$ are recalled in Section 3.1.

The following  example shows that the smash product $A\#H$ is a CY algebra while $A$ itself is not.

\begin{eg}\label{eg1}

Let $H$ be $U(\mc{D},\lmd)$ with the datum $(\mc{D},\lmd)$   given
by \begin{itemize}
\item $\bgm=\lan g\ran$, a free abelian group of rank 1;
\item The Cartan matrix is of type $A_1\times A_1$;
\item $g_1=g_2=g$;
\item $\chi_1(g)=q^2$, $\chi_2(g)=q^{-2}$, where $q$ is not a root of unity;
\item $\lmd_{12}=\frac{1}{q-q^{-1}}$.\end{itemize}
The algebra $H$ is isomorphic to the quantum enveloping algebra $U_q(\mathfrak{sl}_2)$.

 Let $A=\kk\lan u,v\ran/(uv-qvu)$ be the quantum plane.
There is an  $H$-action on $A$ as follows:
$$\begin{array}{cccc}
x_1\cdot u=0,&x_2\cdot u=qv,&g\cdot u=qu,\\
x_1\cdot v=u,&x_2\cdot v=0,&g\cdot v=q^{-1}v.
\end{array}$$
The algebra $A\#H$ is isomorphic to the quantized symplectic oscillator algebra of rank 1 \cite{gakh}.

It is well known that the algebra $A$ is a twisted CY algebra with Nakayama automorphism $\mu$ given by
$$\mu(u)=qu,\;\;\mu(v)=q^{-1}v,$$ and the algebra $H$ is a CY Hopf algebra (\cite[Theorem 3.3.2]{c}). One can easily check that the homological determinant of the $H$-action is trivial and for any $x\in A\#H$, $[\mu\#S^{-2}](x)=gxg^{-1}$. That is, the automorphism $\mu\#S^{-2}$ is an inner automorphism. Therefore, $A\#H$ is a CY algebra.

The invertible elements of $A\#H$ are $\{g^m\}_{m\in \ZZ}$. Therefore, one can see that the automorphism $\id\#S^2$ of $A\#H$ can not be an inner automorphism, although, $S^2$ is an inner automorphism of $H$.

\end{eg}
More generally, we have the following example.
\begin{eg}
Let $H$ be $U(\mc{D},\lmd)$ with the datum $(\mc{D},\lmd)$ given by
\begin{itemize}
\item $\bgm=\lan y_1,y_2,\cdots, y_n\ran$, a free abelian group of rank $n$;
\item The Cartan matrix $\mathbb{A}$ is of type $A_n\times A_n$;
\item $g_i=g_{n+i}=y_i$, $1\se i\se n$;
\item $\chi_i(g_j)=q^{a_{ij}}$, $\chi_{n+i}(g_j)=q^{-a_{ij}}$, $1\se i\se n$, where $q$ is not a root of unity;
\item $\lmd_{ij}=\dt_{n+i,j}\frac{1}{q-q^{-1}}$, $1\se i< j\se 2n$.\end{itemize}
Then $H$ is isomorphic to the algebra $U_q(\mathfrak{sl}_n)$. It is also a CY Hopf algebra.

Let $A$ be the quantum polynomial algebra $$\kk\lan u_1,u_2,\cdots,u_{n+1}\;|\;u_ju_i-qu_iu_j,1\se i<j\se n+1\ran.$$

There is an  $H$-action on $A$ as follows:
$$\begin{array}{ll}x_i\cdot u_j=\dt_{ij}u_{i+1}, 1\se i\se n;&x_i\cdot u_j=\dt_{i+1,j}qu_{i}, n+1\se i\se 2n\\
y_i\cdot u_j=\begin{cases}q^{-1}u_j,&j=i;\\qx_j,&j=i+1;\\x_j,&otherwise.\end{cases}\end{array}$$

It is well known that the algebra $A$ is a twisted CY algebra with Nakayama automorphism $\mu$ given by $\mu(u_i)=q^{n+2-2i}u_i$, $1\se i\se n+1$.

One can also check that the homological determinant of the $H$-action is trivial. The automorphism $\mu\#S^{-2}$ is an inner automorphism. For any $x\in A\#H$, $[\mu\#S^{-2}](x)=gxg^{-1}$, where $g=y_1^ny_2^{2n-2}\cdots y_i^{in-i(i-1)}\cdots y_n^{n^2-n(n-1)}$. Therefore, $A\#H$ is a CY algebra.

Let $H^0$ be the algebra $U(\mc{D},0)$. The algebra $H$ is a cocycle deformation of $U(\mc{D},0)$. Actually, $H\cong (H^0)^\sg$, where $\sg$ is a 2-cocycle on $H^0$ such that $\sg(h_1,h_2)=1$, $\sg(x_i,h_1)=\sg(h_2,x_i)=0$, for all $h_1$, $h_2\in \bgm$ and $1\se i\se n+1$, and
$$\sg (x_i,x_j)=\begin{cases}\lmd_{ij},&j=n+i;\\0,&otherwise.\end{cases}$$
Then we have the crossed product $A\#_\sg H^0$. By Theorem \ref{main}, $A\#_\sg H^0$ is a twisted CY algebra with Nakayama automorphism  $\eta$ defined by $\eta(a\#h)=\mu(a)\#h$, for all $a\#h\in A\#H$. In fact, $\eta$ is an inner automorphism. For any $x\in A\#_\sg H^0$, $\eta(x)=gxg^{-1}$. So $A\#_\sg H^0$ is also a CY algebra.
\end{eg}

\begin{eg}

 Let $H=U(\mc{D},\lmd)$, where $(\mc{D},\lmd)$  is the datum  given
by
\begin{itemize}
\item $\bgm=\lan y_1,y_2\ran$, a free abelian group of rank 2;
\item The Cartan matrix $\mathbb{A}$ is of type $A_1\times A_1$;
\item $g_1=y_1,g_2=y_2$;
\item $\chi_1(g_1)=q^2,\chi_1(g_2)=q^{-4},\chi_2(g_1)=q^4,\chi_2(g_2)=q^{-2}$, where $q$ is not a root of unity;
\item $\lmd=0$.\end{itemize}
The algebra $H$ is a twisted CY algebra with homological integral $_{\xi_1}\kk$, where $\xi_1$ is the algebra homomorphism given by
$$\xi_1(g_1)=q^6g_1, \xi_1(g_2)=q^{-6}g_2, \t{ and } \xi_1(x_i)=0, i=1,2.$$

Let $\sg$ be a 2-cocycle  on $H$ such that $\frac{\sg(g_1,g_2)}{\sg(g_2,g_1)}=q^3$,  $\sg(x_i,g_j)=\sg(g_j,x_i)=0$, $1\se i,j\se 2$,  and
$\sg (x_1,x_2)=\frac{1}{q-1}, \sg (x_2,x_1)=0$.  Then the cocycle deformation  $H^\sg$ is just the algebra $U(\mc{D}',\lmd')$, where
$(\mc{D}',\lmd')$  is the datum given
by \begin{itemize}
\item $\bgm=\lan y_1,y_2\ran$, a free abelian group of rank 2;
\item The Cartan matrix is of type $A_1\times A_1$;
\item $g_1=y_1,g_2=y_2$;
\item $\chi_1(g_1)=q^{-2},\chi_1(g_2)=q,\chi_2(g_1)=q^{-1},\chi_2(g_2)=q^{2}$, where $q$ is not a root of unity;
\item $\lmd_{12}=\frac{1}{q-1}$.\end{itemize}

The algebra $H^\sg$ is a twisted CY algebra with homological integral $_{\xi_2}\kk$, where $\xi_2$ is the algebra homomorphism given by
$$\xi_2(g_1)=q^{-3}g_1, \xi_2(g_2)=q^3g_2, \t{ and } \xi_2(x_i)=0, i=1,2.$$

 Let $A=\kk\lan u,v\ran/(uv-q^2vu)$ be the quantum plane.
There is an  $H^\sg$-action on $A$ as follows:
$$\begin{array}{ccccc}
x_1\cdot u=0,&x_2\cdot u=v,&g_1\cdot u=q^{-1}u,&g_2\cdot u=q^2u\\
x_1\cdot v=u,&x_2\cdot v=0,&g_1\cdot v=qv,&g_2\cdot v=q^{-2}v.
\end{array}$$

We have mentioned in Example \ref{eg1} that $A$ is a twisted CY algebra with Nakayama automorphism $\mu$ given by
$$\mu(u)=q^2u, \;\;\mu(v)=q^{-2}v.$$ One can check that the homological determinant of the $H$ action is trivial. Now we can form the algebras $A\#H^\sg$ and $A\#_\sg H$. By Theorem \ref{cysma}, the algebra $A\#H^\sg$ is a twisted CY algebra with Nakayama automorphism $\mu\#(S^{-2}\circ [\xi]^r)$. This automorphism cannot be an inner automorphism. That is,  $A\#H^\sg$ is not a CY algebra. Theorem \ref{main} shows that the algebra the algebra $A\#_\sg H$ is a twisted CY algebra with Nakayama automorphism $\rho$ defined by $\rho(a)=\mu(a)$, $a\in A$, $\rho(x_1)=q^{-2}x_1$, $\rho(x_2)=q^{2}x_2$, and $\rho(g_i)=\xi(g_i)g_i$, $i=1,2$. The automorphism $\rho$ is an inner automorphism. For any $x\in A\#_\sg H$, $\rho(x)=(g_1^2g_2^2)^{-1}x(g_1^2g_2^2)$. Therefore, the algebra $A\#_\sg H$ is a CY algebra.

\end{eg}

\subsection*{Acknowledgement} This work is supported by an FWO-grant and grants from NSFC (No. 11301126), ZJNSF (No. LQ12A01028).

\vspace{5mm}

\bibliography{}

\end{document}